\documentclass[10pt]{article}
\newcommand{\Bilder}{y}  % include figures with 'y'
\usepackage[colorlinks=true, linktocpage=true, linkcolor=red!70!black, citecolor=green!50!black, urlcolor=blue!50!black]{hyperref}
\usepackage{verbatim}
\usepackage[multiple]{footmisc}
\usepackage{amsfonts,color}
\usepackage{longtable}
\usepackage{makecell, multirow, tabularx}
\usepackage{amsmath,amssymb}
\usepackage{adjustbox}
\usepackage{epsfig}
\usepackage{lscape}
\usepackage{amssymb}
\usepackage{graphicx}
\usepackage[utf8]{inputenc}	
\usepackage{amsthm,amssymb}
\usepackage{aligned-overset}
\usepackage[english]{babel}
\usepackage{dsfont}
\usepackage{bbm}
\usepackage[T1]{fontenc}
\usepackage{colonequals}
\usepackage{geometry}
\usepackage{mathtools}
\usepackage{csquotes}
\usepackage{url}
\usepackage{lipsum}
\usepackage{array}
\usepackage{tablefootnote}
\usepackage{wrapfig}
\usepackage[utf8]{inputenc}
\usepackage[dvipsnames]{xcolor}
\usepackage{accents}
\usepackage{marvosym}
\usepackage{centernot}
\usepackage{tikz}
\usepackage{titling}
\usepackage{multicol}
\newcommand{\blue}[1]{{\color{blue}#1}}

\newcommand{\red}[1]{{\color{red}#1}}

\renewcommand{\phi}{\varphi}
\allowdisplaybreaks[1]
\makeindex
\newtheorem{Theorem}{Theorem}[section]

\theoremstyle{definition}
\newtheorem{Example}[Theorem]{Example}
\newtheorem{Remark}[Theorem]{Remark}

\renewcommand{\rho}{\varrho}

\newcommand{\dxi}{\,\mathrm{d}\xi}

\newcommand{\ddt}{\frac{\mathrm{d}}{\mathrm{d}t}}

\newcommand{\pat}{\partial_t}

\newcommand{\nn}{\nonumber}
\newcommand{\io}{\int_{\Omega_\xi}}

\newcommand{\OO}{{\rm O}}

\newcommand{\ee}{\mathrm{e}}
\newcommand{\Lin}{\mathrm{Lin}}

\newcommand{\R}{\mathbb{R}}

\newcommand{\tr}{\mathrm{tr}}
\newcommand{\GL}{\mathrm{GL}}
\newcommand{\SYM}{\text{\rm Sym}(3)}

\newcommand{\Div}{\mathrm{Div}}%for matrices
\renewcommand{\div}{\mathrm{div}}%for vectors
\newcommand{\sym}{\mathrm{sym}}
\newcommand{\skw}{\mathrm{skew}}
\newcommand{\Id}{{\boldsymbol{\mathbbm{1}}}}
\newcommand{\id}{{\boldsymbol{\mathbbm{1}}}}
\newcommand{\subsectionid}{{\text{\boldmath$\mathbbm{1}$}}}
\newcommand*\dif{\mathop{}\!\mathrm{d}}
\newcommand{\norm}[1]{\lVert #1 \rVert}
\def\dvg{\textnormal{Div}}
\def\H{\mathbb{H}}
\def\C{\mathbb{C}}
\def\sk{\textnormal{skew}}
\def\dd{\displaystyle}
\DeclareMathOperator{\ZJ}{ZJ}
\DeclareMathOperator{\GN}{GN}
\DeclareMathOperator{\Sym}{Sym}
\DeclareMathOperator{\iso}{iso}
\DeclareMathOperator{\Cof}{Cof}
\DeclareMathOperator{\NH}{NH}
\DeclareMathOperator{\lin}{lin}
\DeclareMathOperator{\diag}{diag}
\DeclareMathOperator{\dev}{dev}
\newcommand{\DD}{\mathrm{D}}
\newcommand{\WW}{\mathrm{W}}

\thanksmarkseries{arabic}

\makeatletter
\DeclareFontFamily{OMX}{MnSymbolE}{}
\DeclareSymbolFont{MnLargeSymbols}{OMX}{MnSymbolE}{m}{n}
\SetSymbolFont{MnLargeSymbols}{bold}{OMX}{MnSymbolE}{b}{n}
\DeclareFontShape{OMX}{MnSymbolE}{m}{n}{
    <-6>  MnSymbolE5
   <6-7>  MnSymbolE6
   <7-8>  MnSymbolE7
   <8-9>  MnSymbolE8
   <9-10> MnSymbolE9
  <10-12> MnSymbolE10
  <12->   MnSymbolE12
}{}
\DeclareFontShape{OMX}{MnSymbolE}{b}{n}{
    <-6>  MnSymbolE-Bold5
   <6-7>  MnSymbolE-Bold6
   <7-8>  MnSymbolE-Bold7
   <8-9>  MnSymbolE-Bold8
   <9-10> MnSymbolE-Bold9
  <10-12> MnSymbolE-Bold10
  <12->   MnSymbolE-Bold12
}{}

\let\llangle\@undefined
\let\rrangle\@undefined
\DeclareMathDelimiter{\llangle}{\mathopen}%
                     {MnLargeSymbols}{'164}{MnLargeSymbols}{'164}
\DeclareMathDelimiter{\rrangle}{\mathclose}%
                     {MnLargeSymbols}{'171}{MnLargeSymbols}{'171}
\makeatother
\begin{document}
\title{Rate-form equilibrium for an isotropic Cauchy-elastic formulation: Part I: modeling}
\author{
Patrizio Neff\thanks{
Patrizio Neff, University of Duisburg-Essen, Head of Chair for Nonlinear
Analysis and Modelling, Faculty of Mathematics, Thea-Leymann-Stra{\ss}e 9,
D-45127 Essen, Germany, email: patrizio.neff@uni-due.de}
, \qquad
Nina J.~Husemann\thanks{
	Nina J. Husemann, University of Duisburg-Essen, Chair for Nonlinear
	Analysis and Modelling,  Faculty of Mathematics, Thea-Leymann-Stra{\ss}e 9,
	D-45127 Essen, Germany, email: nina.husemann@stud.uni-due.de}
, \qquad 
Sebastian Holthausen\thanks{
Sebastian Holthausen, University of Duisburg-Essen, Chair for Nonlinear
Analysis and Modelling,  Faculty of Mathematics, Thea-Leymann-Stra{\ss}e 9,
D-45127 Essen, Germany, email: sebastian.holthausen@uni-due.de}
, \\[0.8em]
Franz Gmeineder\thanks{
	Franz Gmeineder, Department of Mathematics and Statistics, University of Konstanz, Universitätsstrasse 10, 78457 Konstanz,
	Germany, email: franz.gmeineder@uni-konstanz.de}
, \qquad 
Thomas Blesgen\thanks{
Thomas Blesgen, Bingen University of Applied Sciences, Berlinstra{\ss}e 109,
D-55411 Bingen, Germany, email: t.blesgen@th-bingen.de}
\\[0.8em]
%, \qquad and \qquad
%N. N.\thanks{
%some text}
}
\maketitle
\vspace{-0,6cm}
\begin{abstract}
% For a specific Cauchy elastic isotropic, compressible constitutive law suitable
% for large rotations and finite strains, we prove the existence of a local weak
% equilibrium solution in terms of stresses and velocities in the current
% configuration at arbitrary smooth compatible Cauchy pre-stress. This is done by converting
% the problem first into a hypoelastic rate-format and by determining explicitly
% the corresponding induced tangent stiffness tensor $\H^{\ZJ}(\sigma)$ acting
% on the corotational Zaremba-Jaumann rate. The induced stiffness tensor
% $\H^{\ZJ}(\sigma)$ turns out to be positive definite for the choices made.
% Then we apply a parabolic regularisation, use the properties of a differential
% Sylvester equation together with the positive definiteness of
% $\H^{\ZJ}(\sigma)$ and pass to the limit together with
% a Schauder fixed-point argument. 

We derive the rate-form spatial equilibrium system for a nonlinear Cauchy elastic formulation in isotropic finite-strain elasticity. For a given explicit Cauchy stress-strain constitutive equation,  we determine those properties that pertain to the appearing fourth-order stiffness tensor. Notably, we show that this stiffness tensor $\H^{\ZJ}(\sigma)$ acting on the Zaremba-Jaumann stress rate is uniformly positive definite. We suggest a mathematical treatment of the {ensuing} spatial PDE-system which may ultimately lead to a local existence result, to be presented in part II of this work. As a preparatory step, we show existence and uniqueness of a subproblem based on Korn's first inequality and the positive definiteness of this stiffness tensor. The procedure is not confined to Cauchy elasticity, however in the Cauchy elastic case, most theoretical statements can be made explicit.

Our development suggests that looking at the rate-form equations of given
Cauchy-elastic models may provide
additional insight to the modeling of nonlinear isotropic elasticity. This 
especially concerns  constitutive requirements emanating from the
rate-formulation, here being reflected by the positive definiteness of $\H^{\ZJ}(\sigma)$.\\

\medskip

\noindent \textbf{Keywords:} nonlinear elasticity, hyperelasticity,
rate-formulation, Eulerian setting, hypo-elasticity, Cauchy elasticity,
material stability, updated Lagrangean, total Lagrangean

\medskip
\noindent \textbf{Mathscinet} classification
74B20 (Nonlinear elasticity),
35A01 (Existence for PDEs),
%74H20 (Existence of Solutions of dynamical problems in solid mechanics)
%%35A01(Existence for PDEs)
%%74B20(Nonlinear elasticity),
%%74H20(Existence of Solutions of dynamical problems in solid mechanics)}
\end{abstract}

%\begin{keyword}
%Willis model \sep metamaterials \sep existence of weak solutions.
%\end{keyword}

%=================================================================

\clearpage

\tableofcontents

\section{Introduction}
\begin{figure}[t]
\begin{center}
\begin{minipage}[h!]{0.9\linewidth}
\if\Bilder y
\centering
\includegraphics[scale=0.6]{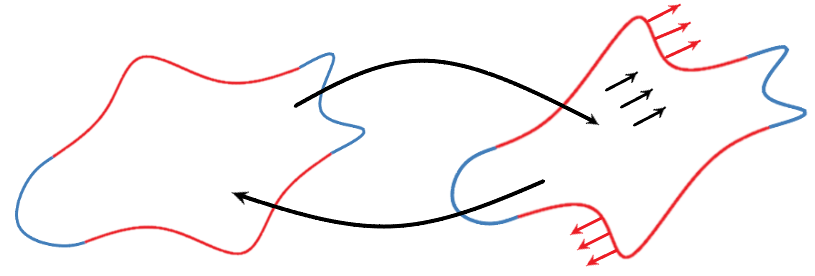}
\put(-300,45){$x \in \Omega$}
\put(-110,45){$\xi \in \Omega_{\xi}$}
\put(-65,75){$f$}
\put(-170,105){$\varphi(x)$}
\put(-210,10){$\varphi^{-1}(\xi)$}
\put(-340,70){\red{$\Gamma_N$}}
\put(-380,22){\blue{$\Gamma_D$}}
\centering
\fi
\caption{Illustration of the domain $\Omega$ with different types of prescribed
boundary conditions and its deformation under the diffeomorphism $\varphi$. The part $\Gamma_D$ describes the Dirichlet boundary, while $\Gamma_N$ defines the Neumann boundary.
%\textcolor{green!50!black}{Isn't this wrong? The red part should be the Neumann part, and the blue one should be the Dirichlet part. \red{$\surd$} }
}
\label{fig8}
\end{minipage}
\end{center}
\end{figure}
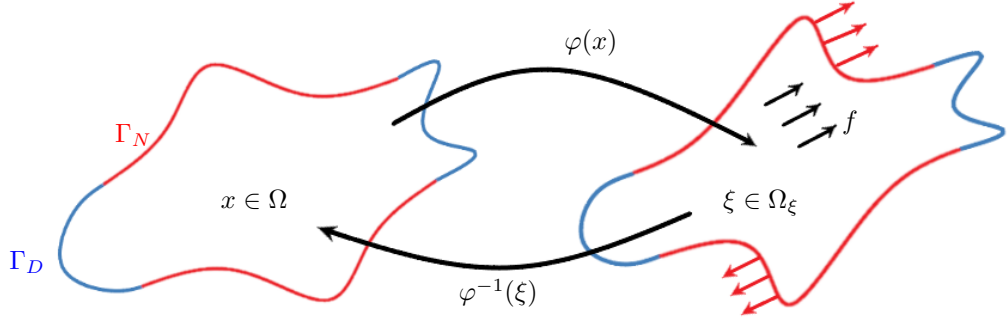

Problems of nonlinear elastic deformations have attracted a great deal of
attention from the applied mathematics community, above all due to their diverse and deep
challenges for analysis. In this work, we focus on homogeneous
isotropic compressible finite strain elasticity. The modelling framework
of these problems can be considered to be complete, see e.g. the classical
books by Ogden \cite{Ogden83}, Ciarlet \cite{ciarlet2022} and
Marsden-Hughes \cite{Marsden83}. Indeed, starting with a hyperelastic
formulation, the static equilibrium equations appear as Euler-Lagrange
equations of the energy functional
\begin{align}
\label{eqintro01}
\int_{\Omega} \WW(\DD\phi) - \widetilde{f} \cdot \phi \dif \text{x} \quad
\longrightarrow \quad \min. \, \phi.
\end{align}
Here, $\phi: \Omega \subset \R^3 \to \R^3$ is the (diffeomorphic)
deformation, $F = \DD \phi (x) \in \GL^+(3)$ is the deformation gradient, and $\widetilde{f}$ represents a dead load
body force in the referential setting. In particular, the elastic energy $\WW: \text{GL}^+(3)\to\R$ completely describes the
constitutive response. Finding suitable forms of $\WW$ which encode physically
realistic behaviour of actual materials is
 \emph{Truesdell's Hauptproblem} \cite{truesdell1956ungeloste}, which, to the best of our knowledge, is not
completely solved up to date. The nonlinear Euler-Lagrange equations
corresponding to \eqref{eqintro01} are then given by 
\begin{equation}
\label{eqintro02}
\begin{alignedat}{2}
\Div_x \, S_1(F(x)) &= \widetilde{f}(x) \qquad && \text{in} \quad \Omega, \\
S_1(F(x)).n &= \widetilde{S}_1(x).n \, = \widetilde{p}(x) \qquad && \text{on} \quad \Gamma_N, \\
\phi(x) &= \widetilde{\phi}(x) \qquad && \text{on} \quad \Gamma_D \, , \quad \text{with } \Gamma = \partial \Omega = \Gamma_D \cup \Gamma_N \, .
\end{alignedat}
\end{equation}
%\blue{$\Gamma = \partial \Omega = \Gamma_D \cup \Gamma_N$}
In \eqref{eqintro02},
%\eqref{eqintro002}
$\Gamma_D \subset \partial \Omega$ is the part of the boundary where
Dirichlet conditions are prescribed,  and $\Gamma_N$ is the part of the boundary
where Neumann conditions apply (see  Figure \ref{fig8}). Here,
$S_1(F) = \DD_FW(F)$ is the non-symmetric first Piola-Kirchhoff stress tensor and $\widetilde{S}_1$ or $\widetilde{p}$ allows to specify the tractions. 

Based on the Piola transformation, it is possible to obtain the equilibrium
equations in the current configuration $\Omega_\xi = \varphi(\Omega)$
by introducing the Cauchy stress tensor
\begin{align}
\sigma(F) = \frac{1}{J} \, S_1(F) \cdot F^T = \frac{1}{J} \, \DD_F \WW(F)
\cdot F^T, \qquad J = \det F, \qquad F = \DD \varphi.
\end{align}
The equilibrium equations in the Eulerian setting read
for $\sigma(\xi) := \sigma(F(\varphi^{-1}(\xi)))$ by abuse of notation
\begin{equation}
\label{eqintro002}
\Div_\xi \, \sigma(\xi) = f(\xi) \qquad \text{in} \quad \Omega_\xi, 
\end{equation}
which can be obtained from \eqref{eqintro02} by setting $f(\xi) = (\det F(x))^{-1} \cdot \widetilde{f}(x)$. We note that this comes to the effect of losing the variational structure of \eqref{eqintro02}, a discrepancy which would not be seen by linearized elasticity.

Objectivity (frame-indifference) and isotropy together amount to
\begin{align}
\WW(Q_1 \, F \, Q_2) = \WW(F), \qquad \forall \, Q_1, Q_2 \in \OO(3), 
\end{align}
and imply that the Cauchy stress $\sigma$ is symmetric and
admits the representation
\begin{align}
\sigma = \sigma(B) = \frac{2}{J} \, \DD_B \widetilde{\WW}(B) \cdot B, \qquad 
B := F \, F^T, \qquad \WW(F) := \widetilde{\WW}(B),
\end{align}
where $\sigma$ is an isotropic tensor function of the left Cauchy-Green tensor $B$, i.e.
\begin{align}
\sigma(Q^T \, B \, Q) = Q^T \, \sigma(B) \, Q, \qquad \forall \, Q \in \OO(3).
\end{align}
For future reference, we note that any such $\sigma$ can be written as (Richter's representation \cite{graban2019,
neff2020, richter1948isotrope, richter1949hauptaufsatze, Richter50, Richter52},
Rivlin-Ericksen representation \cite{BakerEri54})
\begin{align}
\label{eq1.8}
\sigma(B) = \beta_0 \, \id + \beta_1 \, B + \beta_{-1} \, B^{-1},
\end{align}
where $\beta_0, \beta_1, \beta_{-1}$ are scalar-valued functions of the
principal invariants $I_1, I_2, I_3$ of $B$, with
\begin{equation}
\begin{alignedat}{2}
    I_1(B) &= \tr (B) = \norm{F}^2 = \lambda_1^2 + \lambda_2^2 + \lambda_3^2 \, ,\\
    I_2(B) &= \tr (\Cof B) = \norm{\Cof F}^2 = \lambda_1^2 \lambda_2^2 + \lambda_2^2 \lambda_3^2 + \lambda_3^2 \lambda_1^2 \, , \\
    I_3(B) &= \det B = (\det F)^2 = \lambda_1^2 \cdot \lambda_2^2 \cdot \lambda_3^2.
\end{alignedat}
\end{equation}
Here, $\lambda_{1}^{2},\lambda_{2}^{2},\lambda_{3}^{2}$ are the 
%singular values 
eigenvalues of $B$. Based on these considerations, the objective of the present paper is to model isotropic, compressible nonlinear elasticity in a rate format in an Eulerian configuration, see \eqref{eqcompleteintro} below. In particular, this specific formulation gives us access to the positive definiteness of a constitutive fourth-order tangent stiffness tensor $\mathbb{H}^{\mathrm{ZJ}}(\sigma)$. In order to put our modeling into its natural context, we pause to discuss related approaches to nonlinear elasticity first.

\subsection{Mathematical approaches towards nonlinear elasticity}
\subsubsection{The direct method of the calculus of variations}
%
% \blue{(Klar ist, dass es ziemlich verschiedene indirekte Methoden der VR gibt. Ich habe z.B. den Inhalt des Buches von Dacorogna so verstanden, dass es nur eine direkte Methode gibt, aber der Nachweis der Eigenschaften mit verschiedenen Methoden und Annahmen nachgewiesen werden kann.)}\\
It is by now standard to use the framework of the calculus of variations
to show the existence of energy minimizers to \eqref{eqintro01} via the direct method. Ball's seminal introduction of polyconvexity \cite{Ball77}  (together with growth conditions and coercivity estimates)  is the fundamental notion here.
Polyconvexity means  that the elastic energy $\WW: {\rm GL}^+(3) \to \R$
can be written as
\begin{align}
\WW(F) = \mathcal{P}(F, \Cof F, \det F),
\end{align}
where $\mathcal{P}: \R^{3 \times 3} \times \R^{3 \times 3} \times \R \to \R$
is convex and $\Cof F$ is defined by $F^T \, \Cof F = \det F \cdot \id$.
Polyconvexity is sufficient for weak lower semicontinuity in spaces of weakly differentiable functions, 
and implies rank-one convexity of the energy, i.e., 
\begin{align}\label{eq:LegendreHadamardIntro}
\mathrm{D}_F^2 \WW(F).(\xi \otimes \eta, \xi \otimes \eta) \ge 0,
\qquad \forall \, \xi, \eta \in \R^3,
\end{align}
if $\WW \in C^2({\rm GL}^+(3);\,\R)$.
Alternatively, one may express this as rank-one monotonicity of the first Piola-Kirchhoff stress $S_1$, in the sense that
\begin{equation}
	\label{eq:rank_one_S_1}
	\langle S_1(F + \xi \otimes \eta) - S_1(F) , \xi \otimes \eta \rangle \geq 0 \quad \forall \xi, \eta \in \R^3 \setminus \{0\} \, .
\end{equation}
Dunn \cite{Dunn1983} has reformulated this requirement in terms of the Cauchy stress.
Condition \eqref{eq:LegendreHadamardIntro} and \eqref{eq:rank_one_S_1} are also referred to as \emph{Legendre-Hadamard ellipticity} (LH-ellipticity for brevity). A notable condition between polyconvexity and rank-one convexity is \emph{Morrey's quasiconvexity}, meaning that 
	\begin{align}
	\int_{\Omega} \WW(F+ \DD \vartheta) \dif \text{x} \ge \int_{\Omega} \WW(F) \dif \text{x} \qquad \forall \, \vartheta \in C_0^{\infty}(\Omega;\R^{3})\;\forall\,F\in\R^{3\times 3}. 
	\end{align}
This means that the homogeneous configuration is energy-optimal with respect to perturbations that leave affine-linear  boundary data $x \mapsto \varphi(x) = F.x + b$ invariant.  
%\textcolor{red}{Quasiconvexity in itself does not support an existence result under physically relevant conditions.} 
%\textcolor{green!50!black}{We should elaborate on this more; as it stands, it sounds cumbersome. Maybe something in the following direction:}
Even though quasiconvexity is often equivalent to lower semicontinuity of numerous energies on spaces of weakly differentiable functions, it is often too general to support existence results under physically relevant conditions. Being based on the minors of matrix, the main benefit of polyconvexity is that it gives us direct access to the physical requirement that extreme stretches ($\det F \to 0^+$) should be
accompanied by extreme elastic energy ($\WW(F) \to +\infty$), so that energy
minimizers satisfy $\det \DD \varphi(x) > 0$ almost everywhere. \\

\begin{figure}[t]
\begin{center}
\begin{minipage}[h!]{0.4\linewidth}
\if\Bilder y
\centering
\includegraphics[scale=0.45]{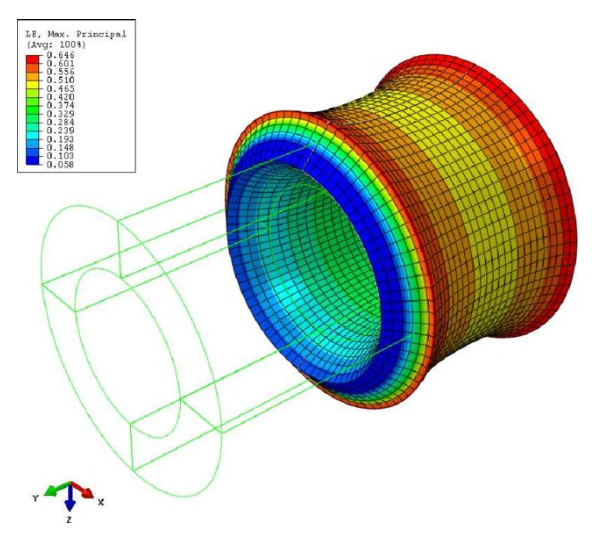}
\caption{The fully everted elastic tube with maximum occurring principal
stretches from \cite{nedjar2018}, finite element solution based on
$\WW(F) \sim \ee^{\norm{\log V}^2}$, $V = \sqrt{F \, F^T}$, the exponentiated Hencky energy.}
\label{fig:everted_elastic_tube}
%\caption{The fully everted elastic tube with maximum occurring principal
%stretches, from \cite{nedjar2018}, FEM-solution based on \break
%$W(F) \sim \ee^{\norm{\log V}^2}$.}
\end{minipage}
\qquad
\begin{minipage}[h!]{0.4\linewidth}
\includegraphics[scale=0.75]{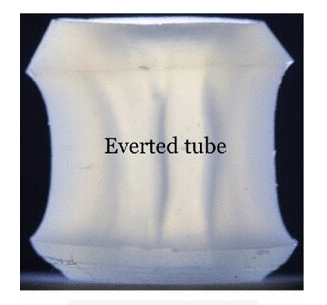}
\fi
\caption{Photo of an everted tube, cf. \cite{nedjar2018}. The everted tube is
not a global energy minimizer, since for traction-free boundary conditions 
the natural state has less energy, cf. Truesdell \cite[p. 103]{truesdell1966}, but
presumably it is a stable local energy minimizer.}
\label{fig:photo_everted_tube}
\end{minipage}
\end{center}
\end{figure}

\noindent The general drawback of the method is that, at present, no general result is available that would show
that the equilibrium equations \eqref{eqintro02} are satisfied even weakly. This is due to the fact that one cannot, in general,  find minimizers satisfying $\det \DD \varphi(x) \ge c^+ > 0$ in advance, 
not to speak of the regularity of the solution. From the perspective of modeling, another issue
is that certain stable equilibria found in nature are not global minimizers. For instance, this is the case for the everted
configuration of a tube, cf. Nedjar et al. \cite{nedjar2018}; see Figures \ref{fig:everted_elastic_tube} and \ref{fig:photo_everted_tube}.
\subsubsection{Implicit function theorem}
A different approach towards existence of solutions to \eqref{eqintro02} is based on the implicit function theorem in appropriately
chosen Sobolev spaces (cf.~Stoppelli \cite{stoppelli1954}, Marsden-Hughes \cite[p.~371]{Marsden83},
Ciarlet \cite{ciarlet2022, ciarlet2010} and Valent \cite{valent1987}).
This methodology establishes the existence of smooth $W^{2,p}$-solutions, $p>3$, to \eqref{eqintro02}
in the neighbourhood of a given smooth equilibrium configuration
$\varphi_0\in C^2(\Omega;\, \R^3)$ with $\det \DD \varphi_0(x) \ge c^+ > 0$ under two fundamental conditions:  
\begin{itemize}
\item uniform Legendre-Hadamard ellipticity at the given configuration $\varphi_0$:
\begin{align}
\label{eqintro001}
\mathrm{D}_F^2 \WW(\DD \varphi_0(x)).(\xi \otimes \eta, \xi \otimes \eta)
\ge c^+ \, |\xi|^2 \, |\eta|^2,
\qquad \forall \, \xi, \eta \in \R^3 \setminus \{0\}, 
\end{align}
\item and the linearization at $\varphi_0$ must be well-posed, i.e.,
\begin{align}
\label{eqintro03}
\Div \, \mathrm{D}_F^2 \WW(\DD \varphi_0(x)).\DD u = 0, \qquad
u \in H_0^1(\Omega, \R^3)
\end{align}
has only the trivial solution. 
\end{itemize} Condition \eqref{eqintro03} excludes any (interior)
instability, and is usually satisfied only in the neighbourhood of
the identity. On the other hand,  \eqref{eqintro001} effectively excludes shear-band type
instabilities as an occurrence of microstructure. In general, this method applies only to either pure Dirichlet problems or pure traction problems.
\subsubsection{Topological continuation methods}
Thirdly, there is the topological continuation method (cf.~Healey
\cite{Healey2000, Healey1997, Healey1998}) which uses degree-theoretical
arguments to show the existence of global $C^{2,\alpha}$-Hölder-smooth
equilibrium solutions to the pure Dirichlet problem \eqref{eqintro02}
driven by body forces parametrized by $t \ge 0$. In this approach, one considers that the elastic energy
$\WW,$ the load $f$ and the domain $\Omega$ are sufficiently smooth, that $\WW$ is strongly
elliptic in $F$, i.e.,~$\mathrm{D}_F^2 \WW(F).(\xi \otimes \eta, \xi \otimes \eta)
> 0$ for all unit vectors
$\xi, \eta \in \R^3 \setminus \{0\}$, that the reference configuration
is stress-free and has a positive definite elasticity tensor $\mathbb{A}:= \DD_F^2 \WW(\id)$.
Then the pure Dirichlet problem in the reference configuration
\begin{align}
\Div_x \, \DD_F \WW(\DD \varphi(x)) = f(x, \varphi, \DD \varphi, t),
\qquad \qquad \quad \varphi(x) = x \qquad \text{on} \qquad \partial \Omega
\end{align}
admits a \emph{solution in the large} (in the language of Healey \cite{Healey2000, Healey1997, Healey1998}) in spaces of smooth functions, meaning that
either the solution exists for all load parameters $t \in \R$, or the
solution explodes at some definite load parameter. In this sense, Healey calls
it a \emph{global inverse function theorem}. \\
\\
%\blue{Am Anfang einmal die LH-Elliptizität ausschreiben, also Legendre-Hadamard-ellipticity (LH-ellipticity). Danach ist die Abkürzung klar.}
Note that in all described existence
theorems some form of LH-ellipticity is crucially involved, allowing to
treat large classes of elastic energies. Altogether, not even for a compressible Neo-Hooke model do we know the
existence of weak equilibrium solutions in situations of interest with mixed
boundary conditions. Here, to the contrary, the existence theory that we prepare will be largely independent of
LH-ellipticity, but will work only for very specific, yet natural examples of constitutive
requirements. These, however, are compatible with the assumptions of nonlinear elasticity for extreme deformations, which usually turn out problematic in a variety of other theories. 
%\textcolor{green!50!black}{but we are getting ahead of ourselves}.

%
\subsection{Isotropic Cauchy-elasticity} \label{isocauchyelast}
Truesdell has extended and generalized the hyperelastic framework by giving up
the requirement that a strain energy $\WW$ exists. While the physics of a theory
without strain energy is certainly questionable, nevertheless, insight into the
nonlinear elasticity problem can be gained and useful results for purely
mechanical problems may be shown \cite{yavari2024}. Here, one directly postulates a
Cauchy stress-stretch constitutive law $B \mapsto \sigma(B)$ and with
\begin{align}
S_1(F) = \sigma(B) \cdot \Cof F, \qquad \qquad B = F \, F^T,
\end{align}
one can either study the problem in the Lagrangian setting \eqref{eqintro02}
or in the Eulerian setting \eqref{eqintro002}. Constitutive assumptions on
$\sigma$ are now called upon for a decent modeling framework. Richter
\cite{graban2019, neff2020, richter1948isotrope, richter1949hauptaufsatze,
Richter50, Richter52} and Reiner \cite{reiner1948elasticity} have early suggested that
\begin{align}
\label{eqintro04}
\sigma: \Sym^{++}(3) \to \Sym(3), \qquad B \mapsto \sigma(B),
\end{align}
should be bijective for idealized perfect elasticity, in a straight-forward extension of
ideas from linear elasticity in which $\sigma_{\lin}: \Sym(3) \to \Sym(3), \;
\varepsilon \mapsto \sigma_{\lin}(\varepsilon) = \C^{\iso} . \varepsilon$
is bijective for $\varepsilon = \sym \, \DD u$. This is, e.g., the case for a
slightly compressible Neo-Hooke type solid with elastic energy
\begin{equation}
\label{eqNeoHooke01}
\begin{alignedat}{2}
\WW_{\NH}(F) &= \frac{\mu}{2} \, \left(\frac{\norm{F}^2}{(\det F)^{\frac23}}
- 3\right) + \frac{\kappa}{2} \, \mathrm{e}^{(\log \det F)^2}, \\
\sigma_{\NH}(B) &= \mu \, (\det B)^{-\frac56} \, \dev_3 B + \frac{\kappa}{2} \,
(\det B)^{-\frac12} \, (\log \det B) \, \mathrm{e}^{\frac14 \,
(\log \det B)^2} \cdot \id,
\end{alignedat}
\end{equation}
where $\dev_n X := X - \frac{1}{n} \, \tr(X) \cdot \id$ (cf.~\cite{CSP2024}). However, most
well-known elastic energies do not satisfy bijectivity of \eqref{eqintro04}
for all positive definite stretches $B \in \Sym^{++}(3)$, e.g.~the slightly simple compressible Neo-Hooke model
\begin{equation}
\begin{alignedat}{2}
\widetilde \WW_{\NH}(F) &= \frac{\mu}{2} \, \left(\frac{\norm{F}^2}{(\det F)^{\frac23}}
- 3\right) + \frac{\kappa}{2} \, (\det F - 1)^2, \qquad \widetilde \sigma_{\NH}(B) &= \mu \, (\det B)^{-\frac56} \, \dev_3 B + \kappa \, (\sqrt{\det B} - 1) \cdot \id
\end{alignedat}
\end{equation}
does not yield bijectivity, due to the merely quadratic nature of the volumetric term \textcolor{red}{\cite{korobeynikov2025}}.

Another frequently-used physical requirement are the (weak) empirical inequalities
(cf. Mihai and Goriely \cite{angela2013numerical, mihai2017,mihai2011positive}).
If $\sigma(B)$ is given as in \eqref{eq1.8}, so
\begin{align}
\sigma(B) = \beta_0 \cdot \id + \beta_1 \cdot B + \beta_{-1} \cdot B^{-1},
\end{align}
it is required that $\beta_1 > 0, \; \beta_{-1} \le 0$ (cf. \cite{thiel2019}). This seems to be in accordance with most available experimental data. Notably, a
positive Poynting effect (cf.~Zurlo et al.~\cite{zurlo2020}) is predicted,
i.e.\;a tube in torsion will lengthen. The Baker-Ericksen inequalities
\cite{BakerEri54} also make a statement for the Cauchy stress (but only at one given
configuration). However, they consider the principal Cauchy stresses and
demand that the higher principal Cauchy stress is occurring for the
higher principal stretch, i.e.
\begin{align}
\label{eqintro06}
(\sigma_i - \sigma_j) \cdot (\lambda_i - \lambda_j) > 0 \quad \iff \quad (\widehat{\sigma}_i - \widehat{\sigma}_j) \cdot (\log \lambda_i - \log \lambda_j) > 0 \, .
\end{align}
In this paper, we will also deal with a version of \, ``stress increases with
strain'',  which we interpret firstly as a condition on the Cauchy-stress
tensor $\sigma$ (the ``true'' stress) and secondly as a subsequent monotonicity
requirement in the Frobenius scalar product (\textbf{Hilbert-monotonicity} of $\widehat \sigma(\log B) \colonequals \sigma(B)$ in terms of $\log B$), cf.~\cite{CSP2024}, where $\log B = 2 \, \log V$ is the logarithmic strain (Hencky strain, true-strain)
\begin{align}
\label{eqintro08}
\langle \widehat{\sigma}(\log B_1) - \widehat{\sigma}(\log B_2), \log B_1
- \log B_2 \rangle > 0, \qquad \forall \, B_1, B_2 \in \Sym^{++}(3),
\quad B_1 \neq B_2,
\end{align}
which is not equivalent to Hilbert-monotonicity of $\sigma$ in terms of $B$, i.e.
\begin{align}
\label{eqintro07}
\langle \sigma(B_1) - \sigma(B_2), B_1 - B_2 \rangle > 0, \qquad \forall \,
B_1, B_2 \in \Sym^{++}(3), \quad B_1 \neq B_2.
\end{align}
For more information about the logarithm of a matrix we refer to
Richter \cite{Richter50}, Lankeit et al.~\cite{LankeitNeffNakatsukasa} and
Neff et al. \cite{Neff_Nagatsukasa_logpolar13}. However,
\eqref{eqintro08} already implies the Baker-Ericksen inequalities and bijectivity of
$B \mapsto \sigma(B)$ if for $\widehat{\sigma}: \Sym(3) \to \Sym(3)$,
$\widehat{\sigma}(\log B):= \sigma(B)$ we have
\begin{equation}
\norm{\widehat{\sigma}(\log B)} \to + \infty \qquad \text{for} \qquad
\norm{\log B} \to + \infty.
\end{equation}
This will be shown in a forthcoming work \cite{MartinVossGhibaNeff}. A
strengthening of the condition \eqref{eqintro08} is the TSTS-M$^{++}$ (true-stress true-strain monotonicity condition)
\begin{equation}
\label{eqCPSdef}
\begin{alignedat}{2}
%\forall \,D\in\Sym(3),\,D\neq0:
% ~
%\langle \frac{\DD^{\ZJ}}{\DD t}[\sigma], D \rangle> 0
%\quad
%&\iff
%\quad 
&\log B \longmapsto \widehat\sigma(\log B)
\;\textrm{is strongly Hilbert-monotone} \\
\iff \quad &\sym \, \DD_{\log B} \widehat \sigma(\log B) \in \Sym^{++}_4(6),
\end{alignedat}
\end{equation}
where $\Sym^{++}_4(6)$ denotes the set of positive definite, minor and major symmetric fourth\textcolor{red}{-}order tensors.\\
Recently, in \cite{Leblond2024} it is shown that TSTS-M$^{++}$ is equivalent to the \textbf{c}orotational \textbf{s}tability \textbf{p}ostulate (CSP)
\begin{equation}
    \big\langle \frac{\DD^{\circ}}{\DD t}[\sigma] , D \big\rangle > 0 \quad \forall \ D \in \Sym(3) \setminus \{0\} \, ,
\end{equation}
see also \cite{poscor2024,secondorderwork2024}, where $\frac{\DD^{\circ}}{\DD t}[\sigma] = \frac{\DD}{\DD t} [\sigma] + \sigma \, \Omega - \Omega \, \sigma$ with $\Omega \in \mathfrak{so}(3)$ denotes any corotational derivative of the Cauchy stress $\sigma$. The corotational stability postulate implies that in all uniform deformations without rotation of principal axes the incremental Cauchy stress moduli are positive \cite{secondorderwork2024}.
%\textcolor{green!50!black}{Explain what $[\sigma]$ is!}
%
\subsection{Hypoelasticity: rate-formulation in the current configuration}
Finally, Truesdell extended and formalized the concept of Cauchy elasticity
towards a rate-formulation. This means the constitutive law is written in terms
of an ODE that must be integrated along a specific loading path for each
material point. Such an ODE must satisfy the principle of frame-indifference
(objectivity) and necessitates the introduction of objective time rates for
stresses. One such stress rate (among many others) is the Zaremba-Jaumann rate 
(cf. Jaumann \cite{jaumann1905, jaumann1911geschlossenes}, Zaremba \cite{zaremba1903forme} and many other
authors \cite{fiala2009, fiala2011, fiala2016, fiala2020objective, kolev2024objective,
romano2011, romano2014}). For its derivation, we now interpret the
diffeomorphism $\varphi$ as a time-dependent quantity, so that
$\varphi = \varphi(x,t): \R^3 \times [0,\infty) \to \R^3$.
Then the corotational Zaremba-Jaumann rate is given by
\begin{align}
\label{eqintroZJrate}
\frac{\DD^{\ZJ}}{\DD t}[\sigma] \colonequals  \frac{\DD}{\DD t}[\sigma] + \sigma \, W - W \,
\sigma, \qquad W = \skw \, L, \qquad L = \dot{F} \, F^{-1} \, ,
\end{align}
where $\frac{\DD}{\DD t}[\sigma]$ denotes the material derivative of $\sigma$.
%\textcolor{green!50!black}{What is $[\sigma]$?}
$\frac{\DD^{\ZJ}}{\DD t}$ is arguably the simplest objective rate since the spin $W$ can be determined
directly from $\DD_\xi v = L = D + W, \, W = \skw \, \DD_\xi v$,
where $v$ is the spatial velocity in the current configuration $\Omega_\xi$.
The Zaremba-Jaumann rate of the Cauchy stress $\frac{\DD^{\ZJ}}{\DD t}[\sigma]$ provides a measure
of the change of stress in a material point with respect to a (corotated) frame
that rotates with the spin $W$, such that rigid body rotations are properly taken
into account. Then a hypoelastic model in the Eulerian configuration
$\Omega_\xi$ can be written as
\begin{equation}
\label{eqintro09}
\begin{alignedat}{2}
\underbrace{\frac{\DD^{\ZJ}}{\DD t}[\sigma] = \H^{\ZJ}(\sigma).D, \qquad D = \sym \,
\DD_\xi v}_{\text{constitutive law}}, \qquad \quad \underbrace{\Div_\xi \, \sigma(\xi,t) = f(\xi, t)}_{\text{spatial equilibrium}}, \qquad \text{in} \quad \Omega_\xi,
\end{alignedat}
\end{equation}
together with suitable initial and boundary conditions. Note that
$\H^{\ZJ}(\sigma)$, defined by \eqref{eqintro09}$_1$, is a constitutive fourth-order tangent stiffness tensor,
mapping symmetric arguments to symmetric arguments (minor symmetry) and depending on the used stress rate, here the Zaremba-Jaumann rate.\\

We will show subsequently that for a smooth diffeomorphism $\varphi \colon \Omega \to \Omega_{\xi}$ the system \eqref{eqintro09} is equivalent to\\
%\fbox{
\begin{minipage}[t][1.8cm][c]{\linewidth}
\begin{equation}
\label{eqcompleteintro}
\boxed{
\begin{alignedat}{2}
\pat \sigma + \DD_\xi \sigma.v + \sigma \, W - W \, \sigma
&= \H^{\ZJ}(\sigma).\sym \, \DD_\xi v, \qquad \quad W = \sk \, \DD_\xi v, \\
\Div_{\xi}[\H^{\ZJ}(\sigma).\sym \, \DD_\xi v] &= \Div_\xi \left[\sigma \cdot
(\DD_\xi v)^T - \div_\xi \, v \cdot \sigma + \sigma \, (\sk \, \DD_\xi v)
- (\sk \, \DD_\xi v) \, \sigma\right] \\
&\qquad \qquad + \div_\xi \, v \cdot f(\xi,t) + \DD_{\xi} f(\xi,t).v
+\pat f(\xi,t),
\end{alignedat}
}
\end{equation}
\end{minipage}
%}

\vspace{0.3cm}
\noindent where only the stress rate and the spatial velocity appear. 

\medskip
\noindent It is well-known that any hyperelastic or Cauchy elastic model
in which $\sigma: \Sym^{++}(3) \to \Sym(3), \break B \mapsto \sigma(B)$ is
bijective can be written in the format \eqref{eqintro09}
(cf.~Truesdell \cite{truesdellremarks} and Noll \cite{Noll55}, see also   {\cite{CSP2024}}), however with an
expression for $\H^{\ZJ}(\sigma)$ that is not easily manageable (cf.~\cite{CSP2024}).
Hypoelasticity of grade zero, which means $\H^{\ZJ}(\sigma) = \C^{\iso} =$ const. (cf. Truesdell \cite{truesdell55hypo}) is usually not
compatible with Cauchy elasticity (the Cauchy stress would depend on the
loading path and not only on the local configuration), and not compatible with
the existence of a strain energy function. This is why
hypoelasticity (also referred to as hypoelasticity of grade zero) has {often been}  abandoned in favour of the total Lagrangean approach in
\eqref{eqintro01}, \eqref{eqintro02}. Yet, as we argue below, \eqref{eqcompleteintro} arises from natural modeling requirements and comes with a satisfactory existence theory.\\ 

\noindent 
For more details on hypoelasticity, we refer the reader to 
\cite{Bernstein60, Green56, korobeynikov2024, Ortiz83,Renardy90, romano2011,
romano2014} for a non-exhaustive list. Moreover, different aspects of the FEM-implementation of hypoelastic models
are addressed e.g.~in \cite{bellini2015, federico2022, govindjee1997,
korobeynikov2023, Ortiz83, Zohdi2006}, and the recent work \cite{CSP2024} provides more in-depth explanations of hypoelasticity.
\subsection{Approach in this paper} \label{sec1.4}
Based on our above discussion, the main purpose of the present paper is to understand whether writing the equations of isotropic,
compressible nonlinear elasticity in the rate format \eqref{eqcompleteintro}
%\red{diese equation Nummer stimmt, das ist die Gleichung um die es geht, sollte in der Intro evtl nochmal aufgenommen werden} 
and the
Eulerian configuration can support any new existence theorem, provided that suitable
constitutive assumptions on the induced tangent stiffness tensor $\H^{\ZJ}(\sigma)$ are
made. Therefore, we first need to reformulate \eqref{eqintro09} into the equivalent system \eqref{eqcompleteintro}
that allows us to utilize the symmetry and positive definiteness of $\H^{\ZJ}(\sigma)$, see Section \ref{alternatesystem}.

In order to be sufficiently self-contained, we will confine ourselves to
one exceptionally {favorable} Cauchy-elastic material for which most necessary
calculations and observations can be made explicit. Namely, we consider the isotropic and compressible Cauchy
elastic constitutive law
\begin{align}
\label{eqintro10}
\sigma: \Sym^{++}(3) \to \Sym(3), \qquad \sigma(B) = \frac{\mu}{2} \,
(B- B^{-1}) + \frac{\lambda}{2} \, \log \det B \cdot \id = \mu \, \sinh (\log B) + \frac{\lambda}{2} \, \tr(\log B) \, \id \, , 
\end{align}
in which $\frac12 (B - B^{-1})$ is the Mooney-strain (cf.~Curnier and
Rakotomanana \cite{curnier1991}). Thus, our development
should only be taken as a first step for more general considerations, including
the most important case of hyperelasticity.

\begin{Remark}
\label{rem1}
The material response \eqref{eqintro10}, while definitely not hyperelastic in the compressible
case \cite{yavari2024b}, presents a number of salient mechanical features, listed below in {Section \ref{rem3}}. Notably, \eqref{eqintro10} is valid for large rotations
(since frame-indifferent) and finite strains (since extreme stretches go with
extreme Cauchy stresses). Moreover, the corresponding induced tangent stiffness tensor
$\H^{\ZJ}(\sigma)$ is smooth as a function of $\sigma$ and not only minor symmetric, meaning that it maps the space of symmetric matrices onto itself, but also major symmetric, whereby $\H^{\ZJ}(\sigma)$ is even self-adjoint; see the appendix, Section \ref{sec:minmajsym} for more detail. Lastly, $\H^{\ZJ}(\sigma)$ is uniformly positive definite, {even though} 
\eqref{eqintro10} is not Legendre-Hadamard elliptic throughout. In the sense that the rank-one monotonicity \eqref{eq:rank_one_S_1} is not satisfied throughout. Note that \eqref{eq:rank_one_S_1} remains applicable to merely Cauchy elastic response \cite{Dunn1983}. Moreover, we note that the law \eqref{eqintro10} may already be used for the stress analysis in a purely mechanical context.
\end{Remark}

While the present paper is centered around the underlying mathematical modeling, the follow-up paper \cite{blesgen2024} will address the corresponding existence theory. In particular, the aim of \cite{blesgen2024} is to establish that the new system \eqref{eqcompleteintro} -- to be derived in Section \ref{alternatesystem} -- 
together with the material response \eqref{eqintro10} and
suitable initial and boundary conditions admits a weak solution in terms of
Cauchy stresses $\sigma$ and the velocity $v$. For the complete initial boundary value problem \eqref{eqcompleteintro} 
%\blue{\eqref{M1}?} 
in the velocity $v(\xi,t)$ and
the Cauchy stress $\sigma(\xi,t)$, we shall use a parabolic regularisation and pass
to the limit in conjunction with a Schauder fixed point argument. This determines the
velocity field $v(\xi,t)$. In an additional integration step, the deformation
solution $\varphi(x,t)$ can be obtained for small $0<t<T^*$, where $T^*$ is
determined by a Schauder-type fixed point argument provided that the velocity solution $v$ is sufficiently
smooth. Yet, in Theorem \ref{lem1}, below we foreshadow such results by establishing the well-posedness of an associated subproblem.\\ 

\noindent 
For technical simplicity, the present and the follow-up paper \cite{blesgen2024} will primarily be concerned with a fixed domain
$\Omega_\xi = \mathrm{const.}$, meaning that no shape or volumetric changes at the boundary are considered. We thereby give the foundation for more general scenarios, where the recently developed theory of time-varying domains will play a pivotal role. With the latter previously being confined to problems from fluid mechanics, see e.g. \cite{breit2018} for applications to fluid-structure-interaction, this will eventually lead to a unified existence theory for nonlinear solid and time-varying materials. The case of time-varying domains, however, comes with a multitude of technical challenges on top of the conceptual key novelties,  which are already at the heart of the scenario as considered here.

\section{The non-linear hypoelasticity model} \label{secmodel}
%
%A parabolic regularization is treated in \cite{CK16}.
%
As described in Section \ref{sec1.4}, the first major task consists in the derivation of a system that is equivalent to \eqref{eqintro09} and allows for the utilization of the symmetry and positive definiteness of the induced fourth-order tangent stiffness tensor $\H^{\ZJ}(\sigma)$. Therefore, let us consider $\WW(F)$ to be the strain energy function of an elastic
material in which $F(x) := \DD \varphi(x) \in {\rm GL}^+(3)$ is the gradient
of a deformation $\varphi: \Omega \to \mathbb{R}^3$ with an open domain $\Omega\subset\mathbb{R}^3$.
The deformation $\varphi(x)$ maps from a stress-free reference configuration $\Omega$ to a
configuration $\varphi(\Omega) =: \Omega_{\xi}$ in the Euclidean 3-space; $\WW(F)$ is measured
per unit volume of the reference configuration. 

In the quasi-static setting, which is the setting we are working in, we understand
$\varphi : \Omega \times \R_+^0 \to \varphi(\Omega,t) \subset \R^3$
as a deformation that allows for interior movements depending on the time $t$ due to loads.
We generally assume that for any fixed time $t_0 \ge 0$ the function
$\varphi(\cdot,t_0): \Omega \to \varphi(\Omega,t_0), \; \varphi(x,t_0) = \xi$
is a diffeomorphism from $\R^3 \to \R^3$. Furthermore, we restrict our analysis to the special case
$\varphi(\Omega, t) \equiv \Omega_\xi \equiv\mathrm{const}$. This corresponds to the situation where only
interior movements of the material are possible due to body forces $f(\xi,t)$, similar to
a moving fluid confined in a container. This setting excludes structural and
geometrical instabilities like e.g.~buckling, barrelling and necking, while
material instabilities are not a priori neglected.

\subsection{Derivation of the hypo-elastic system - a rate-formulation in the spatial setting} \label{alternatesystem}

Now, to start with the derivation of the alternative system of equations it is mandatory to obtain an expression for $\Div_\xi \frac{\DD}{\DD t}[\sigma]$, since our goal is to replace $\frac{\DD}{\DD t}[\sigma]$ by
	\begin{align}
	\H^{\ZJ}(\sigma).D = \frac{\DD}{\DD t}[\sigma] + \sigma \, W - W \, \sigma \qquad \iff \qquad \frac{\DD}{\DD t}[\sigma] = \H^{\ZJ}(\sigma).D - \sigma \, W + W \, \sigma \, ,
	\end{align}
yielding an equation with the structure
	\begin{align}
	\Div_{\xi}[\H^{\ZJ}(\sigma(\xi,t)).D] = - g(\xi,t)
	\end{align}
with some vector-valued function $g(\xi,t)$. {For future reference, we note that it is {incorrect}
to simply write $\frac{\DD}{\DD t} [\Div_{\xi}\sigma] = \Div_{\xi}\frac{\DD}{\DD t}[\sigma]$. This will be shown by the subsequent calculations, and is due to the fact that 
$\Div_\xi = \Div_{\varphi(x,t)}$ is itself time-dependent and $\frac{\DD}{\DD t}$
in the spatial configuration is meant to be the material or substantial
derivative}. Thus, we first recall how the corresponding equilibrium
equation
\begin{align}
\label{eqcurrentref01}
\Div_\xi \, \sigma(\xi,t) = f(\xi,t),
\end{align}
where $f(\xi,t)$ is a vector-valued body force, can be transformed from the
current configuration $\Omega_\xi$ to the reference configuration $\Omega$
and vice versa. This will allow us to differentiate the equation  {with respect to} time in the
reference configuration. A key ingredient in this respect is the
\textbf{Piola transformation}, which states that for a diffeomorphism
$\varphi: \Omega \to \Omega_{\xi}$ and a continuously differentiable operator
$S: \Omega_{\xi} \to \R^{3 \times 3}$ the following identity holds: 
\begin{align}
\label{eqpiola1}
\Div_x (S(\varphi(x,t),t) \cdot \Cof F(x,t)) = \det F(x,t) \cdot \Div_{\xi}
S(\varphi(x,t),t) = \det F(x,t) \cdot \Div_{\xi} S(\xi,t).
\end{align}
\subsubsection{Transformation of the equilibrium equation} \label{sec1.1}
Consider the equilibrium equation in elastostatics
	\begin{align}
	\label{eq001}
	\Div_\xi \, \sigma(\xi,t) = f(\xi,t)
	\end{align}
in the current configuration $\Omega_\xi$ for a given body force $f(\xi,t)$. The non-symmetric first Piola-Kirchhoff stress tensor $S_1(x):= S_1(F(x))$ connects points $x \in \Omega$ with points $\varphi(x,t) = \xi \in \Omega_\xi$ via
	\begin{align}
	S_1(x,t) = \sigma(\xi,t) \cdot \Cof F(x,t).
	\end{align}
 Applying the Piola transformation \eqref{eqpiola1} to $\sigma(\xi,t)$ yields
	\begin{equation}
	\label{eq003}
	\begin{alignedat}{2}
	\Div_x \, S_1(x,t) &= \Div_x \big(\sigma(\varphi(x,t),t) \cdot \Cof F(x,t)\big) \overset{\text{Piola}}{=} \det F(x,t) \cdot \Div_\xi \, \sigma(\xi,t) \\
	\iff \qquad \Div_\xi \, \sigma(\xi,t) &= (\det F(x,t))^{-1} \cdot \Div_x \, S_1(x,t).
	\end{alignedat}
	\end{equation}
Simultaneously, we can transform the body force $f = f(\xi,t)$, suppressing the
$t$-dependence for this calculation, via (cf.~Ciarlet \cite[p.73]{ciarlet2022})
\begin{align}
\io \langle \Div_\xi \, \sigma(\xi), \vartheta(\xi) \rangle \dif \xi &=
\io \langle f(\xi), \vartheta(\xi) \rangle \dif \xi \qquad \forall \, \vartheta
\in C^\infty_0(\Omega_{\xi}, \R^3) \notag \\
&= \int_{\Omega} \langle f(\varphi(x)), \vartheta(\varphi(x)) \rangle \cdot
\underbrace{\det \DD \varphi(x) \dif \text{x}}_{\cong \,\dif \xi} = \int_{\Omega}
\langle \underbrace{\det \DD \varphi(x) \, f(\varphi(x))}_{=: \,
\widetilde{f}(x)}, \underbrace{\vartheta(\varphi(x))}_{=: \,
\widetilde{\vartheta}(x)} \rangle \dif \text{x} \nn \\
\label{eq002}
&= \int_{\Omega} \langle \widetilde{f}(x), \widetilde{\vartheta}(x) \rangle
\dif \text{x} \qquad \forall \widetilde{\vartheta} \in C^{\infty}_0 (\Omega, \R^3),
\end{align}
whereby we have, in particular, the relation 
\begin{align}
\label{eq004}
\widetilde{f}(x,t) = \det \DD \varphi(x,t) \cdot f(\varphi(x,t)) =
\det F(x,t) \cdot f(\xi,t),\;\;\;\text{so}\;\;\; f(\xi,t) = (\det F(x,t))^{-1}
\cdot \widetilde{f}(x,t).
\end{align}
Thus, using \eqref{eq003} and \eqref{eq004}, we can express \eqref{eq001}
equivalently in the reference configuration $\Omega$ by
\begin{align}
\label{eq01}
\Div_x \, S_1(F(x,t)) = \widetilde{f}(x,t) = \det F(x,t) \cdot f(\xi,t).
\end{align}
\subsubsection{Derivation of the rate formulation of force equilibrium for
the Zaremba-Jaumann rate}
\label{sec1.2}
Now that we have an equivalent system to \eqref{eqcurrentref01} in the
reference configuration $\Omega$, we can interchange the operations as
initially intended, i.e.~$\ddt \Div_x S_1 = \Div_x \ddt S_1$, since in the
material setting the domain $\Omega$ is independent of the time $t$. Even though clear to the expert reader, the following derivation of the rate formulation is crucial for the sequel and thus shall be carried out in detail: Since the first Piola-Kirchhoff stress tensor $S_1(x,t)$
fulfills the relation
\begin{align}
S_1(x,t) = \sigma(\varphi(x,t),t) \cdot \Cof F(x,t) = J(x,t) \cdot
\sigma(\xi,t) \cdot F^{-T}(x,t),
\end{align}
where $J = \det F$, we have
\begin{align}
\label{eq02}
\ddt S_1 = \ddt [\det F \cdot \sigma \cdot F^{-T}] = \ddt [\det F] \cdot \sigma
\cdot F^{-T} + \det F \cdot \frac{\DD}{\DD t} [\sigma]\cdot F^{-T} + \det F \cdot \sigma
\cdot \ddt [F^{-T}] \, .
\end{align}
Additionally, we obtain
\begin{align}
\ddt [\det F] = \langle \Cof F, \dot F \rangle = \det F \, \langle F^{-T},
\dot F \rangle = \det F \, \langle \id, \dot F \, F^{-1} \rangle = \det F
\cdot \tr(L)
\end{align}
as well as
\begin{align}
\ddt [F^{-T}] = \left(\ddt F^{-1}\right)^T = (- F^{-1} \, \dot F \, F^{-1})^T
= - L^T \, F^{-T}.
\end{align}
Hence, \eqref{eq02} yields the relation
\begin{equation}
\label{eq03}
\begin{alignedat}{2}
\ddt S_1 &= \det F \cdot \tr(L) \cdot \sigma \cdot F^{-T} + \det F \cdot
\frac{\DD}{\DD t} [\sigma] \cdot F^{-T} + \det F \cdot \sigma \cdot (-L^T \cdot F^{-T}) \\
&= \tr(D) \cdot \sigma \cdot \Cof F + \frac{\DD}{\DD t} [\sigma] \cdot \Cof F - \sigma
\cdot L^T \cdot \Cof F \\
&= \left(\tr(D) \cdot \sigma + \frac{\DD}{\DD t} [\sigma] - \sigma \cdot L^T\right) \cdot
\Cof F.
\end{alignedat}
\end{equation}
Applying the Piola transformation \eqref{eqpiola1} once again to transform back
to the current configuration $\Omega_\xi$, we can use \eqref{eq01} and
\eqref{eq03} to obtain 
\begin{equation}
\label{eq04}
\begin{alignedat}{2}
\Div_x \ddt S_1(F(x,t)) \overset{\eqref{eq03}}&{=} \Div_x \left(\left[\tr(D)
\cdot \sigma + \frac{\DD}{\DD t} [\sigma] - \sigma \cdot L^T \right] \cdot \Cof F \right) \\
\overset{\text{Piola}}&{=} (\det F) \cdot \Div_\xi \left(\tr(D) \cdot \sigma
+ \frac{\DD}{\DD t} [\sigma] - \sigma \cdot L^T \right) = \ddt \widetilde{f}(x,t) \\
\overset{\substack{\det F > 0 \\ \left. \right.}}{\iff} \qquad \Div_{\xi}
\frac{\DD}{\DD t} [\sigma] &= (\det F)^{-1} \cdot \ddt \widetilde{f}(x,t) - \Div_\xi
\left(\tr(D) \cdot \sigma - \sigma \cdot L^T \right),
\end{alignedat}
\end{equation}
where now by \eqref{eq004}
\begin{equation}
\begin{alignedat}{2}
\ddt \widetilde{f}(x,t) &= \langle \Cof F, \dot F \rangle \, f(\varphi(x,t),t)
+ \det F \cdot \ddt [f(\varphi(x,t),t)] \\
&= \det F \cdot \bigg( \underbrace{\langle F^{-T},\dot F \rangle}_{= \, \tr(D)}
\cdot f(\varphi(x,t),t) + \ddt [f(\varphi(x,t),t)]\bigg).
\end{alignedat}
\end{equation}
Thus $\eqref{eq04}_3$ becomes
\begin{align}
\label{eq05}
\Div_{\xi} \frac{\DD}{\DD t} [\sigma] = \tr(D) \cdot f(\varphi(x,t),t)+ \ddt [f(\varphi(x,t),t)]
-\Div_\xi \left(\tr(D) \cdot \sigma - \sigma \cdot L^T \right).
\end{align}
In a final step we use the defining relation for the Zaremba-Jaumann derivative
$\frac{\DD^{\ZJ}}{\DD t}[\sigma]$
\begin{align}
\H^{\ZJ}(\sigma).D = \frac{\DD^{\ZJ}}{\DD t}[\sigma] := \frac{\DD}{\DD t} [\sigma]
+ \sigma \, W - W \, \sigma, \qquad W = \sk L, \qquad L = \dot F \, F^{-1}
\end{align}
for an invertible constitutive law $B \mapsto \sigma(B)$ in \eqref{eq05} and
where $\H^{\ZJ}(\sigma)$ is the induced fourth-order tangent stiffness tensor,
to obtain
\begin{align}
\Div_{\xi}[\H^{\ZJ}(\sigma).D] &= \tr(D) \cdot f(\varphi(x,t),t)
+ \ddt [f(\varphi(x,t),t)] - \underbrace{\Div_\xi \left(\tr(D) \cdot \sigma
-\sigma \cdot L^T \right)}_{\text{``contribution from Piola''}}
+\underbrace{\Div_{\xi}[\sigma \, W - W \,
\sigma]}_{\text{``contribution from the rate''}} \notag \\
\label{eqdefofg}
&= \underbrace{\Div_\xi \left(\sigma \cdot L^T - \tr(D) \cdot \sigma
+\sigma \, W - W \, \sigma\right) + \tr(D) \cdot f(\varphi(x,t),t)
+ \ddt [f(\varphi(x,t),t)]}_{\equalscolon \, - \, g(\xi,t)} \, .
\end{align}
%This means that 
So, we have determined the alternative rate form equilibrium system
in the current configuration
\begin{equation}
\label{eqlong}
\begin{alignedat}{2}
\Div_{\xi}[\H^{\ZJ}(\sigma).\sym \, \DD_\xi v] = \; &\Div_\xi \left[
\sigma \cdot (\DD_\xi v)^T - \div_\xi \, v \cdot \sigma + \sigma \,
(\sk \, \DD_\xi v) - (\sk \, \DD_\xi v) \, \sigma\right] \\
& \qquad + \div_\xi \, v \cdot f(\xi,t) + \DD_{\xi} f(\xi,t).v + \pat f(\xi,t).
\end{alignedat}
\end{equation}
The system is completed with the constitutive equation for the symmetric
Cauchy stress $\sigma$
\begin{align}
\frac{\DD^{\ZJ}}{\DD t}[\sigma] = \H^{\ZJ}(\sigma).\sym \, \DD_\xi v \qquad
\iff \qquad \pat \sigma + \DD_\xi \sigma.v + \sigma \, W - W \, \sigma
= \H^{\ZJ}(\sigma).\sym \, \DD_\xi v,
\end{align}
in which $W = \sk \, \DD_\xi v$ and $\sigma \in \Sym(3)$. \\
\\
Thus, the complete rate-form system in the Eulerian setting for $\sigma(\xi,t)\in\Sym(3)$
and $v(\xi,t) \in \R^3$ is given by \\
%\fbox{
\begin{minipage}[t][1.8cm][c]{\linewidth}
\begin{equation}
\label{eqcomplete}
\boxed{
\begin{alignedat}{2}
\pat \sigma + \DD_\xi \sigma.v + \sigma \, W - W \, \sigma
&= \H^{\ZJ}(\sigma).\sym \, \DD_\xi v, \qquad \quad W = \sk \, \DD_\xi v, \\
\Div_{\xi}[\H^{\ZJ}(\sigma).\sym \, \DD_\xi v] &= \Div_\xi \left[\sigma \cdot
(\DD_\xi v)^T - \div_\xi \, v \cdot \sigma + \sigma \, (\sk \, \DD_\xi v)
- (\sk \, \DD_\xi v) \, \sigma\right] \\
&\qquad \qquad + \div_\xi \, v \cdot f(\xi,t) + \DD_{\xi} f(\xi,t).v
+\pat f(\xi,t)
\end{alignedat}
}
\end{equation}
\end{minipage}
%}
\\[0.8em]
as already announced in \eqref{eqcompleteintro}. In particular, we note that the deformation $\varphi(x,t)$ and therefore $F = \DD \varphi$ and
$B = F \, F^T$ do not appear explicitly anymore. As a consequence, the deformation $\varphi(x,t)$
has to be found in a subsequent step, the latter being potentially conditional upon sufficient
%enough
regularity of the velocity $v(\xi,t)$.
%
%(Neue Notation)
\begin{Remark} \label{remnotation}
This system, paired with suitable initial and boundary conditions, will be the
origin for the existence theory to be developed in the follow-up paper \cite{blesgen2024}.
%Section \ref{sec2.30}. 
At this point, it is
important to take careful note of the following perceived ambiguity:
While it is necessary for the derivation of the system~\eqref{eqcomplete} to
understand $\xi = \varphi(x,t)$ as a \textbf{time-dependent} variable,
leading for example to the relation for the substantial derivative
\begin{align}
\frac{\DD}{\DD t} [f(\xi,t)] = \ddt [f(\varphi(x,t),t)] = \pat f(\xi,t) + \DD_\xi f(\xi,t).
\varphi_t = \pat f(\xi,t) + \DD_\xi f(\xi,t).v \, ,
\end{align}
it is evident that from now on, working only on the current configuration,
$\xi$ is understood as \textbf{time-independent, solely spatial variable}.
As correctly pointed out by an anonymous referee, this approach can be encountered frequently in the context of fluid mechanics. In particular, every time-derivative of $\sigma = \sigma(\xi,t)$ or $v = v(\xi, t)$
is to be understood as partial derivative in time direction. We will clarify
this causality by writing $\pat \sigma(\xi,t)$ instead of
$\ddt [\sigma(\xi,t)]$, even though, strictly speaking, in the context of the
initial boundary value problem \eqref{M1} both are the same. 
\end{Remark}
\begin{Remark}[Comparison to the linear theory]
Considering zero prestress $\sigma \equiv 0$ we obtain
\begin{align}
\label{eqrem01}
\Div_{\xi} [\H^{\ZJ}(0). \sym \, \DD_\xi v(\xi,t)] = \div_\xi \, v
\cdot f(\xi,t) + \DD_{\xi} f(\xi,t).v + \partial_t f(\xi,t),
\end{align}
which is formally \emph{nearly} the rate formulation
\begin{align}
\Div_x \C^{\iso} . \sym \, \DD_x u(x,t) = f(x,t) \qquad
\overset{\partial_t}{\implies} \qquad \Div_x \C^{\iso} . \sym \, \DD_x u_t(x,t)
=\pat f(x,t)
\end{align}
of linear elasticity by identifying reference and deformed configuration
together with $\H^{\ZJ}(0) = \C^{\iso}$. Note that for spatially homogeneous
body force $f$ and small $\tr(D) = \text{div}_\xi\, v$ (nearly incompressible),
the right hand side of \eqref{eqrem01} reduces to $\pat f(\xi,t)$.
\end{Remark}
\begin{Remark}[Comparable results]
Setting $\ddt [\widehat f] := (\det F)^{-1} \cdot \ddt [\widetilde f(x,t)]$, equation \eqref{eq04} is equivalent to
	\begin{align}
	\label{eq06}
	\Div_\xi \left(\frac{\DD}{\DD t}[\sigma] + \tr(L) \, \sigma - \sigma \, L^T\right) - \ddt [\widehat{f}] = 0.
	\end{align}
This equation is similar to that of many other authors. \\
\\
For example Ji et al.~\cite[Eq.~44]{ji2013} write
(in the weak form and in our notation)
	\begin{equation}
            \label{eq:Ji}
		\begin{alignedat}{2}
			%		\int_{\Omega_\xi} \Big\langle \underbrace{\frac{\DD^{\ZJ}}{\DD t}[\sigma]
			%		- 2 \, D \, \sigma + L \, \sigma + \tr(L) \, \sigma}_{=: \, A}
			%		\red{- \frac{\dif}{\dif t}[\widehat{f}]}, \DD_{\xi} \vartheta \Big\rangle \dif \xi
			%		= 0 \qquad \forall \, \vartheta \in C^{\infty}_0(\Omega_\xi, \R^3).
			0 &= \int_{\Omega_\xi} \langle \frac{\DD^{\ZJ}}{\DD t} [\sigma] , \sym \, \DD_{\xi} \vartheta \rangle - 2 \langle \sigma , D \, \sym \, \DD_{\xi} \vartheta \rangle + \langle \sigma , L^T \, \DD_{\xi} \vartheta \rangle + \tr(L) \, \langle \sigma , \sym \, \DD_{\xi} \vartheta \rangle \dif \xi \qquad \forall \, \vartheta \in C^{\infty}_0(\Omega_\xi, \R^3).\\
			&= \int_{\Omega_\xi} \langle \frac{\DD^{\ZJ}}{\DD t} [\sigma] - 2 \, D \, \sigma + \tr(D) \, \sigma , \sym \, \DD_{\xi} \vartheta \rangle + \langle L \, \sigma , \DD_{\xi} \vartheta \rangle \dif \xi \\
			&= \int_{\Omega_\xi} \langle \frac{\DD^{\ZJ}}{\DD t} [\sigma] - 2 \cdot \frac{1}{2} \, (D \, \sigma + \sigma \, D) + \tr(D) \, \sigma + L \, \sigma , \DD_{\xi} \vartheta \rangle \dif \xi\\
			&= \int_{\Omega_\xi} \langle \underbrace{\frac{\DD^{\ZJ}}{\DD t} [\sigma] - 2 \cdot \frac{1}{2} \, (D \, \sigma + \sigma \, D) + L \, \sigma + \tr(D) \, \sigma}_{\equalscolon \mathcal{A}} , \DD_{\xi} \vartheta \rangle \dif \xi \qquad \forall \, \vartheta \in C^{\infty}_0(\Omega_\xi, \R^3)\, .
		\end{alignedat}
	\end{equation}
	Rewriting $\mathcal{A}$ then leads to ($\tr(L) = \tr(D)$)
	\begin{equation}
	\label{eq06a}
	\begin{alignedat}{2}
	\mathcal{A} &= \frac{\DD}{\DD t}[\sigma] + \sigma \, W - W \, \sigma - 2 \cdot \frac12 \cdot (D \, \sigma + \sigma \, D) + D \, \sigma + W \, \sigma + \tr(D) \, \sigma \\
	&= \frac{\DD}{\DD t}[\sigma] + \sigma \, W - \sigma \, D + \tr(D) \, \sigma = \frac{\DD}{\DD t}[\sigma] + \tr(D) \, \sigma - \sigma \, L^T,
	\end{alignedat}
	\end{equation}
which is the argument of the divergence proposed in \eqref{eq06}. \\
\\
Another example is given by Aubram in 
%an upcoming paper
% vorher Prop. 5.3
\cite[Prop 6.3]{aubram2017}, who considers (again in the weak form and in our notation)
\begin{align}
	\label{eq07}
	0 = \int_{\Omega_\xi} \Big\langle \frac{\DD^{\ZJ}}{\DD t}[\sigma] - 2 \, D \, \sigma + \sigma \, \tr(D), \sym \, \DD_\xi \vartheta \rangle + \langle \sigma, L^T \cdot (\DD_\xi \vartheta) \Big\rangle \dif \xi, \qquad \forall \, \vartheta \in C^{\infty}_0(\Omega_\xi,\R^3).
\end{align}
Using $\langle \sigma, L^T \cdot \DD_{\xi} \vartheta \rangle = \langle L \, \sigma, \DD_{\xi} \vartheta \rangle$, it is easily verified that \eqref{eq07} coincides with the previous two equations \eqref{eq06} and \eqref{eq:Ji}. \\
\\
Lastly, Korobeynikov and Larichkin \cite[p.~6]{korobeynikov2024} write (with $\underline{\Div}$ denoting their divergence operator)
	\begin{align}
	\underline{\Div} \left(\frac{\DD}{\DD t}[\sigma] + \sigma \, W - W \, \sigma + \sigma \, \tr(D) + \sigma \, L^T - (D \, \sigma + \sigma \, D) \right) + \varrho \, \frac{\dif}{\dif t}[f] = 0,
	\end{align}
which may be simplified by observing
	\begin{align}
	\sigma \, L^T + \sigma \, W - W \, \sigma - D \, \sigma - \sigma \, D = \sigma \, (D - W) + \sigma \, W - W \, \sigma - D \, \sigma - \sigma \, D = -W \, \sigma - D \, \sigma = - L \, \sigma,
	\end{align}
yielding the identity
%equation
	\begin{align}
	\label{eq08}
	\underline{\Div} \left(\frac{\DD}{\DD t}[\sigma] + \sigma \, \tr(D) - L \, \sigma \right) + \varrho \, \frac{\dif}{\dif t}[f] = 0.
	\end{align}
Since in their notation $\underline{\Div}(X):= \Div_\xi(X^T)$, we recover \eqref{eq06} by observing that
	\begin{align}
	\left(\frac{\DD}{\DD t}[\sigma] + \sigma \, \tr(D) - L \, \sigma \right)^T = \frac{\DD}{\DD t}[\sigma] + \sigma \, \tr(D) - \sigma \, L^T.
	\end{align}
We also note here that the other authors assume that the force on the right hand side of the equation is given as a (spatially constant) body force per unit volume, i.e.~we can write it in the form $\widetilde{f}(t)\coloneqq \varrho \, f(t)$, where $\varrho = \frac{\varrho_0}{\det F}$ and $\varrho_0 $ is the density for $t=0$. Then, starting with the spatial equilibrium equation
	\begin{align}
	\label{01}
	\Div_\xi \, \sigma - \widetilde{f}(t)= 0,
	\end{align}
one may proceed with the Piola transformation, which yields
	\begin{align}
	\Div_\xi \, \sigma(\xi,t) = (\det F(x,t))^{-1} \cdot \Div_x \, S_1(x,t).
	\end{align}
Hence \eqref{01} becomes
	\begin{align}
	\frac{1}{\det F} \cdot \Div_x \, S_1 - \widetilde{f}(t) = 0 \qquad \iff \qquad \Div_x \, S_1 - \varrho_0 \cdot f(t) = 0.
	\end{align}
Now differentiation with respect to time yields
	\begin{align}
	\Div_x \ddt S_1 - \varrho_0 \cdot \dot{f}(t) = 0, \qquad \text{with} \qquad \Div_x \ddt S_1 = (\det F) \cdot \Div_\xi \left(\tr(D) \cdot \sigma + \frac{\DD}{\DD t}[\sigma] - \sigma \cdot L^T \right).
	\end{align}
Thus ($\rho = \frac{\varrho_0}{\det F}$),
	\begin{equation}
	\begin{alignedat}{2}
	(\det F) \cdot \Div_\xi \left(\tr(D) \cdot \sigma + \frac{\DD}{\DD t}[\sigma] - \sigma \cdot L^T \right) - \varrho_0 \cdot \dot{f} &= 0 \\
	\iff \qquad \Div_\xi \left(\tr(D) \cdot \sigma + \frac{\DD}{\DD t}[\sigma] - \sigma \cdot L^T \right) - \varrho \cdot \dot f &= 0. \end{alignedat}
	\end{equation}
The approach in this paper is a more general one, in which $\frac{\dif}{\dif t} [\widehat f]$ admits the explicit expression
	\begin{equation}
	\begin{alignedat}{2}
	\ddt [\widehat f] &= \det F \cdot \left(\tr(D) \cdot f(\varphi(x,t),t) + \ddt [f(\varphi(x,t),t)]\right) \\
	&= \det F \cdot \bigg(\div_\xi \, v \cdot f(\xi,t) + \DD_\xi f(\xi,t).v + \partial_t f(\xi,t)\bigg),
	\end{alignedat}
	\end{equation}
where the force $f$ may also change in spatial direction.
\end{Remark}
\begin{Remark}
If we consider $\H^{\ZJ}(\sigma).D = 2 \, \mu \, D + \lambda \, \tr(D) \, \id$ in \eqref{eqcomplete} then we deal with a classical zero-grade hypo-elastic rate-formulation. To the best of our knowledge, no existence result for \eqref{eqcomplete} is known but it is clear that $\H^{\ZJ}(\sigma).D = 2 \, \mu \, D + \lambda \, \tr(D) \, \id$ may not be integrable towards a Cauchy-elastic formulation \cite{simopister}. If we take $\H^{\ZJ}(\sigma)$ as induced tangent stiffness tensor subordinate to the Cauchy-elastic law $\sigma = \frac{\mu}{2} \, (B-B^{-1}) + \frac{\lambda}{2} \, \log \det B \, \id$ then showing existence to \eqref{eqcomplete} may lead to an existence result for the problem
	\begin{align}
	\Div_{\xi} \sigma(\xi,t) = f(\xi,t), \qquad \sigma = \frac{\mu}{2} \, (B-B^{-1}) + \frac{\lambda}{2} \, \log \det B \, \id
	\end{align}
provided that the initial conditions are compatible and the solutions $(\sigma,v)$ to \eqref{eqcomplete} are smooth enough to reconstruct the deformation $\varphi(x,t)$ along the particle moving in the velocity field. We therefore find it worthwhile to first investigate the rate-form equilibrium system under generic conditions on $\H^{\ZJ}(\sigma)$, postponing the regularity issue to future works.
\end{Remark}

%\clearpage
%
\subsection{{Properties of the constitutive law} $\sigma(B) = \frac{\mu}{2} \, (B - B^{-1}) + \frac{\lambda}{2} \,(\log \det B) \cdot \subsectionid$}
\label{rem3}
The constitutive choice of $\sigma(B)$ as made in \eqref{eqintro10} is not arbitrary. In fact, the 
 constitutive law 
\begin{align*}
\widehat \sigma(\log B) := \sigma(B) = \frac{\mu}{2} \, (B - B^{-1}) + \frac{\lambda}{2} \,
(\log \det B) \cdot \id = \mu \, \sinh (\log B) + \frac{\lambda}{2} \, \tr(\log B) \, \id
\end{align*}
shares the following salient properties:
\begin{itemize}
\item The constitutive law $B \mapsto \sigma(B)$ is objective and isotropic. 
%\blue{Was bedeutet es, dass $\sigma(B)$ `objective` ist? Entspricht dies der physikalischen Bedeutung?}
\item $\log B \mapsto \widehat{\sigma}(\log B) = \sigma(B)$ is strictly
monotone in $\log B$ (it satisfies TSTS-M$^{+}$), meaning that
\begin{align}\label{eq:logstrictlymonotone}
\langle \widehat{\sigma}(\log B_1) - \widehat{\sigma}(\log B_2),
\log B_1 - \log B_2 \rangle > 0\qquad\forall B_1, B_2 \in \Sym^{++}(3),\;\;\;\;B_1 \neq B_2. 
\end{align}
However, we point out that $B \mapsto \sigma(B)$ is not monotone in $B$. 
\item In $\H^{\ZJ}(\sigma) . D:= \frac{\DD^{\ZJ}}{\DD t}[\sigma]$ the fourth-order tensor $\H^{\ZJ}(\sigma)$ is 
uniformly positive definite, meaning that there is a constant $c^+ > 0$ such that
$\langle \H^{\ZJ}(\sigma).D,D \rangle \ge c^+ \, \norm{D}^2$ for all $\sigma \in \Sym(3)$ for all $D \in \Sym(3)$, and has major and minor symmetry (cf.~\cite{CSP2024}),
i.e.~$\H^{\ZJ}(\sigma) \in \Sym^{++}_4(6)$.
\item $\sigma(B)$ fulfills the ``tension-compression symmetry'' $\sigma(B) = -\sigma(B^{-1})$.
\item Extreme stresses for extreme strains:
$\!\left\{\!\!
\begin{array}{l}\norm{\sigma(B)} \to + \infty \;\; \text{for} \;
\det B \to 0 \; \text{and} \; \norm{\sigma(B)} \to +\infty \;\; \text{for} \;
\norm{B} \to \infty, \\ \norm{\widehat{\sigma}(\log B)} \to +\infty \;\;
\text{as} \; \norm{\log B} \to + \infty, \end{array}\right.$ 

Therefore, $\sigma(B)$ is suitable for large rotations and large strains.
\item $B \mapsto \sigma(B), \, V \mapsto \sigma(V), \, \log B \mapsto
\widehat{\sigma}(\log B), \, \log V \mapsto \widehat{\sigma}(\log V)$
are all bijective.
\item There exists a smooth inverse mapping for $B \mapsto \sigma(B)$: 
$\mathcal{F}^{-1}: \Sym(3)\to \Sym^{++}(3), \quad B = \mathcal{F}^{-1}(\sigma(B))$.
\item $\sigma(B)$ induces a rank-one convex formulation in a large
neighbourhood of the stress-free reference configuration.
\item Correct linearization: 
$2 \, \mathrm{D}_B \sigma(B) \big\vert_{B=\id} =
\C^{\iso} \in \Sym^{++}_4(6)$, \qquad
$\C^{\iso}.\varepsilon = 2 \, \mu \, \varepsilon + \lambda \,
\tr(\varepsilon) \cdot \id, \quad \varepsilon = \sym \, \DD u$.
\item The weak-empirical inequalities are satisfied ($\beta_1 = \frac{\mu}{2} > 0$ and
$\beta_{-1} = -\frac{\mu}{2} < 0$); no condition on $\beta_0$.
\item The tension-extension inequality, pressure-compression inequality and
Baker-Ericksen inequalities are satisfied (cf.~the Appendix \ref{secA1}) and we have monotonicity in
uniaxial loading and monotonicity of shear stress in simple shear
(cf.~\cite{CSP2024} and \cite{young2025,secondorderwork2024}) .
\item $B \mapsto \sigma(B)$ is additionally operator-monotone in $B$. By
{\it operator-monotonicity} we mean that \break
$B_1 \prec B_2 \implies \sigma(B_1) \prec \sigma(B_2)$, where $B_1 \prec B_2$
is the Löwner partial ordering \cite{loewner1934}:
\begin{align}
B_1 \prec B_2 \qquad \iff \qquad B_2 - B_1 \in \Sym^+(3).
\end{align}
The operator-monotonicity of $B$ and $B^{-1}$ with respect to $B$ is clear (cf.~Löwner~\cite{loewner1934}),
it remains to show operator monotonicity of
$B \mapsto \log \det B \cdot \id = \tr(\log B) \cdot \id$. Since if
$B_1 \prec B_2$ we have $\log B_1 \prec \log B_2$
(due to the operator-monotonicity of $B \mapsto \log B$, it is clear that
$\log B_2 - \log B_1 \in \Sym^{+}(3)$), then
$\tr(\log B_2 - \log B_1) > 0$, so that
$\tr(\log B_1) \cdot \id \prec \tr(\log B_2) \cdot \id$.
\item There is a Cauchy pseudo-stress potential
$\widetilde{\Psi} = \widetilde{\Psi}(B)$ for the Mooney strain
$\frac12 (B - B^{-1})$:
\begin{align}
\underbrace{B - B^{-1}}_{\text{monotone in} \, B} = \mathrm{D}_B
\underbrace{\left(\frac12 \norm{B}^2 - \log \det B +
3\right)}_{\text{convex in} \, B} = \mathrm{D}_B \widetilde{\Psi}(B)
\end{align}
but $\sigma(B) = \frac{\mu}{2} \, (B-B^{-1}) + \frac{\lambda}{2} \,
\log \det B \cdot \id$ is \textbf{not Hilbert-monotone} in terms of $B$, i.e.
\begin{align}
\langle \sigma(B_1) - \sigma(B_2), B_1 - B_2 \rangle \not> 0 \qquad \forall \,
B_1, B_2 \in \Sym^{++}(3), \quad B_1 \neq B_2.
\end{align}
Crucially, comparing this with 
\eqref{eq:logstrictlymonotone}, we particularly find $\sigma$ to be non-(Hilbert) monotone in $B$, but in $\log B$. Moreover, there does not exist a function $\Psi: \Sym^{++}(3) \to \R$ such that
\begin{align}
\sigma(B) = \DD_B \Psi(B),
\end{align}
while $\H^{\ZJ}(\sigma)$ is major symmetric. Moreover, $\sigma$ is
\textbf{not hyperelastic}, i.e.~there does not exist an isotropic energy
$\WW(F)$ so that
\begin{align}
\sigma(B) = \frac{2}{J} \, \DD_B \WW(B) \cdot B \, .
\end{align}
\end{itemize}

\subsection{The complete spatial PDE-system and some first considerations}
\label{sec2.30}
For a final time $0<T\le\infty$ we define $\Omega_T:=\Omega_\xi\times(0,T)$ and
$\Sigma_T:=\partial\Omega_\xi\times(0,T)$.
Then, based on the system \eqref{eqcomplete}, we study the following
hypo-elasticity initial boundary value problem in terms of symmetric
Cauchy stresses $\sigma$ and velocities $v$: \\
\\
\noindent{\it Find the solution
$(\sigma,v) = (\sigma(\xi,t),v(\xi,t))$ with values in  $\Sym(3) \times \R^3 \cong \R^6
\times \R^3 = \R^9$ of}
\begin{equation}
\label{M1}
\begin{alignedat}{2}
\pat \sigma + \DD_\xi \sigma.v + \sigma \, W - W \, \sigma &=
\H^{\ZJ}(\sigma).\sym \, \DD_\xi v, \qquad \quad W = \sk \, \DD_\xi v \in
\mathfrak{so}(3), && \text{in} \quad \Omega_T, \\
\Div_{\xi}[\H^{\ZJ}(\sigma).\sym \, \DD_\xi v] &= \left\{
\begin{array}{l} \Div_\xi\! \left[
\sigma(\DD_\xi v)^T-(\div_\xi \, v)\sigma
+ \sigma \, (\sk \, \DD_\xi v) - (\sk \, \DD_\xi v) \, \sigma\right] \\ 
\qquad \quad + (\div_\xi \, v)f(\xi,t) +(\DD_{\xi}f(\xi,t))v
+\pat f(\xi,t) ,
\end{array} \right.
	&& \text{in} \quad \Omega_T, \\
v(\xi,0) & = v_0(\xi), \qquad \sigma(\xi,0) = \sigma_0(\xi) \in \Sym(3),
&& \text{in} \quad \Omega_\xi, \\
v(\xi,t) &= 0, && \text{on} \quad \! \Sigma_T. \\
\end{alignedat}
\end{equation}
Counting equations, $\eqref{M1}_1$ yields six and $\eqref{M1}_2$
yields three equations for the independent variables $(\sigma,v)$.\\
Note that \eqref{M1} constitutes a system resembling Differential Algebraic Equations (DAEs), a subject area with a rich literature. However, the first equation also contains spatial derivatives, so that DAE-techniques are difficult to be applied directly. We now discuss this system in more detail and give a well-posedness result for a subproblem in Theorem \ref{lem1}, leaving the full system \eqref{M1} to the follow-up paper \cite{blesgen2024}.
%\cmag
%Die Symmetrie von $\sigma$ sollte aus (\ref{S1}) und Lemma~\ref{lem3} folgen,
%wenn $\HH(\sigma).\sym \DD_\xi v$ und $\sigma_0$ symmetrisch sind.
%Dazu: $\sigma W-W\sigma=2\,\sym(\sigma W)$ ist symmetrisch.
%Der Term $\pat\sigma+\DD_\xi\sigma.v$ ist ein Transportterm und macht die
%Symmetrie nicht kaputt!
%Für die Abschätzungen wären zusätzliche $\sym$-Operatoren rechts kein Problem!
%\cn
%\\
%\red{
%Bemerkung: In jeder Approximation wird $\sigma \in \Sym(3)$ gesucht, das ist
%hart dabei, dann $\sigma = \sigma(\xi,t)$. $\DD_\xi \sigma(\xi,t)$ ist ein
%dreistufiger Tensor $\DD_\xi \sigma(\xi,t).h \overset{?}{=} \sum_{i=1}^3
%\underbrace{\DD_{\xi_i} \sigma(\xi_1, \xi_2, \xi_3, t)}_{\in \Sym(3)}
%\cdot h_i \in \Sym(3)$ $\forall \, h \in \R^3$!
%}
%
\subsubsection{The induced tangent stiffness tensor $\H^{\ZJ}(\sigma)$}\label{sec:stiffness}
Searching for a solution $(\sigma(\xi,t), v(\xi,t))$ of the system \eqref{M1},
one cannot expect that the initially assumed constitutive law
\begin{align}
\label{eqintro003}
\sigma(B)= \frac{\mu}{2}(B-B^{-1})+\frac{\lambda}{2} (\log \det B) \cdot \Id,
\qquad \text{where} \quad B = F \, F^T \quad \text{is the Finger tensor},
\end{align}
is a-priori fulfilled for every time $0<t<T$, notably if some form of
approximation is used. However, the information which constitutive law was used
to derive the system~\eqref{M1} is expressed by the initial condition
$\sigma(\xi,0) = \sigma_0(\xi)$ and more importantly, in an implicit way,
by the induced tangent stiffness tensor $\H^{\ZJ}(\sigma)$.
As shown in the Appendix~\ref{sec:rate-formulations} (see also \cite{CSP2024}
for an alternative method), the constitutive law \eqref{eqintro003} leads to
the induced tangent stiffness tensor
\begin{align}
\label{eqintro004}
\H^{\ZJ}(B). D := \frac{\mu}{2} \, \{B \, D + D \, B + B^{-1} \, D + D \,
B^{-1}\} + \lambda \, \tr(D) \cdot \id.
\end{align}
Furthermore, since the constitutive law $\sigma \colon \Sym^{++}(3) \to \Sym(3)$ \eqref{eqintro003} is invertible, the
inverse function $\mathcal{F}^{-1}(\sigma) = B$ can be used to derive the
explicit dependence $\sigma \mapsto \H^{\ZJ}(\sigma)$, thus encoding which
constitutive law was used initially to obtain the system \eqref{M1}.
\begin{Remark}
When starting from hyperelasticity, $\H^{\ZJ}(\sigma)$
always possesses \emph{minor symmetry}, meaning that
 $\H^{\ZJ}(\sigma) :\SYM\to\Lin(\SYM,\SYM)$, but no \emph{major symmetry},
i.e.~$\H^{\ZJ}(\sigma)$ is not necessarily symmetric as a matrix. However, 
the particular choice \eqref{eqintro003} for the constitutive law
$B \mapsto \sigma(B)$ leads to an induced fourth-order tangent stiffness tensor
$\H^{\ZJ}(\sigma)$ that is \textbf{minor and major symmetric}
as well as \textbf{positive definite} (cf. Section~\ref{sec:rate-formulations}).
\end{Remark}

\noindent It is important to note that the positive definiteness implies the 
existence of a constant $c_0>0$ with 
\begin{equation}
\label{H1}
\boxed{\big\langle \H^{\ZJ}(\sigma). D, D \big\rangle \ge c_0 \,
\norm{D}^2, \quad \qquad \forall \, D\in\SYM,}
\end{equation}
not to be confused with LH-ellipticity of the constitutive law
(as \eqref{eqintro003} is not LH-elliptic in the sense of \eqref{eq:rank_one_S_1}).
Since \eqref{H1} holds independently of the stress level $\sigma$, this
condition resembles nevertheless the strong ellipticity of $\H^{\ZJ}(\sigma)$,
cf.~\cite{ADN64}.

\begin{Remark}
\label{rem5}
Every hyperelastic formulation in the deformed configuration $\Omega_\xi$ can
be written in the format (\ref{M1}), if the constitutive law
$\sigma: \Sym^{++}(3) \to \Sym(3), \; B \mapsto \sigma(B)$ is invertible
(cf.\;Truesdell \cite{truesdellremarks} and Noll \cite{Noll55}).
However, $\H^{\ZJ}(\sigma)$ need not be symmetric or positive definite 
at given $\sigma$ and the invertibility of $B \mapsto \sigma(B)$
is equivalent almost everywhere in $\sigma$ to the invertibility of
$\H^{\ZJ}(\sigma)$ (cf.~\cite{CSP2024}).

For example, a slightly compressible Neo-Hooke type solid with elastic energy
\begin{equation}
\label{eqpolyconvex001}
\begin{alignedat}{2}
\WW_{\NH}(F) &= \frac{\mu}{2} \, \left(\frac{\norm{F}^2}{(\det F)^{\frac23}}
- 3\right) + \frac{\kappa}{2} \, \mathrm{e}^{(\log \det F)^2}, \\
\sigma_{\NH}(B) &= \mu \, (\det B)^{-\frac56} \, \dev_3 B + \frac{\kappa}{2} \,
(\det B)^{-\frac12} \, (\log \det B) \, \mathrm{e}^{\frac14 \, (\log \det B)^2}
\cdot \id,
\end{alignedat}
\end{equation}
admits an invertible constitutive law
$\sigma_{\NH} : \Sym^{++}(3) \to \Sym(3), \; B \mapsto \sigma_{\NH}(B)$
together with $\det \H^{\ZJ}(\sigma_{\NH}) \neq 0$, but it can be
shown that the corresponding induced tangent stiffness tensor
$\H^{\ZJ}(\sigma)$ is not positive definite throughout and $\WW_{\NH}$ does therefore not satisfy the TSTS-M$^{++}$ condition, and therefore does not satisfy inequality  \eqref{eq:logstrictlymonotone}, see \cite{CSP2024} and \cite{korobeynikov2025}.
However, \eqref{eqpolyconvex001} is polyconvex and LH-elliptic
(cf. Hartmann and Neff \cite{Hartmann2002}). 
\end{Remark}
\begin{Remark}
In the general theory of hypo-elasticity, one considers a rate equation of the
form
\begin{align}
\frac{\DD^{\sharp}}{\DD t}[\sigma] = \H^*(\sigma).D \, ,
\end{align}
where $\frac{\DD^{\sharp}}{\DD t}[\sigma]$ is an appropriate objective rate of
the Cauchy stress tensor $\sigma$ and $\H^*(\sigma)$ is a constitutive
fourth-order tangent stiffness tensor (for more information on this topic see
\cite{CSP2024} and \cite{federico2024}). As the notation already suggests, the choices of
$\frac{\DD^{\sharp}}{\DD t}[\sigma]$ and $\H^*(\sigma)$ are a-priori arbitrary
and independent of each other. In this work, it is our choice to use the
tangent stiffness tensor $\H^{\sharp}(\sigma)$ that is induced by
$\frac{\DD^{\sharp}}{\DD t}[\sigma]$, yielding several features discussed
earlier.

However, since the fourth-order tangent stiffness tensor $\H^*(\sigma)$
can be prescribed arbitrarily in general, one commonly used choice is  
\begin{align}
\label{eqintro005}
\H^*(\sigma).D = \C^{\iso}.D = 2 \, \mu \, D + \lambda \, \tr(D)
\cdot\id, \qquad \text{the so-called ``zero grade hypo-elasticity''.}
\end{align}
We again emphasize that for the choice $\H^* = \C^{\iso}$ there is no explicit law
for the Cauchy stress $\sigma$.
However, $\C^{\iso}.D$ from \eqref{eqintro005} would also be positive definite
for $\mu > 0, 3 \, \lambda + 2 \, \mu > 0$, completely independent of
the stress. 
\end{Remark}
\subsubsection{Initial, boundary and compatibility conditions}
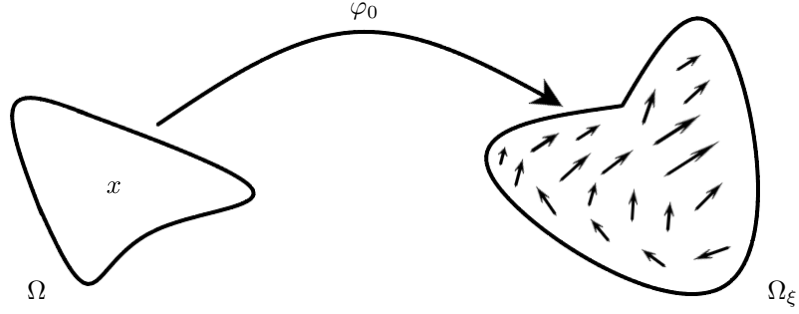
\begin{figure}[t]
\begin{center}
\begin{minipage}[h!]{0.9\linewidth}
\if\Bilder y
\centering
\includegraphics[scale=0.5]{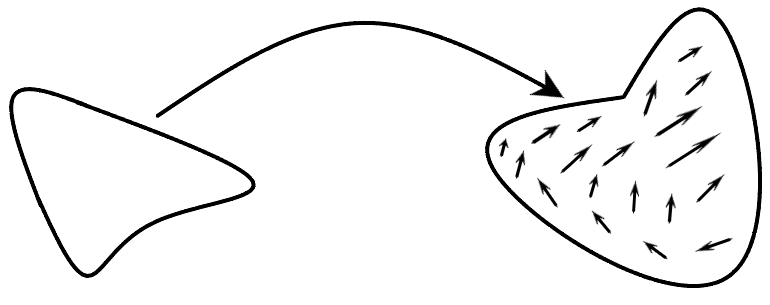}
\put(0,5){$\Omega_{\xi}$}
\put(-280,5){$\Omega$}
\put(-250,45){$x$}
\put(-158,113){$\phi_0$}
\fi
\caption{Picture of the initial diffeomorphism
$\phi_0 = \phi(x,0): \Omega \to \Omega_{\xi}$. The arrows indicate the traction
vectorfield $\sigma_0(\xi). e_1$, where $\sigma_0(\xi)$ is given by
\eqref{eqinitial001}. Additionally, we have $v_0(\xi) = 0$ on
$\partial \Omega_\xi$. 
%\red{$\neq \Sigma_T$?}.
%%denn $\partial \Omega_\xi$ ist der Rand für $t=0$,
%%$\Sigma_T$ ist der Rand für $t \in (0,T)$.
%%\red{Geometrisch(im $\R^3$) ändert sich der Rand von $\Omega_\xi$ aber nicht.}
}
%%\green{Für $t=0$ sollte $\partial\Omega_\xi=\partial\Omega$ gelten ?!?}
%%\red{Nein, das dürfte im allgemeinen nicht gelten. Für $t=0$ deformiert
%%$\varphi_0(x)$ das Gebiet $\Omega$ auf das neue Gebiet $\Omega_\xi$, welches
%%nicht denselben Rand haben muss. Was aber gilt ist, dass der Rand sich
%%im Laufe der Zeit nicht verformt}.
%%\green{Im allgemeinen gilt also $\phi_0=\Id$ nicht ?!}}
\label{fig7}
\end{minipage}
\end{center}
\end{figure}

The equations $\eqref{M1}_3$ and $\eqref{M1}_4$ state initial and boundary
conditions that have to be fulfilled by the solution $(\sigma(\xi,t),v(\xi,t))$
of the system $\eqref{M1}$, where $\sigma_0$ is a (smooth) initial Cauchy
stress distribution and $v_0$ is the initial velocity, respectively.
For an illustration of the initial configuration see Figure~\ref{fig7}.

\noindent The functions $\sigma_0$ and $v_0$ need to fulfill a set of requirements necessary for the application of the Schauder fixed point theory in the follow-up \cite{blesgen2024} and the determination of the deformation $\varphi(x,t)$, which we will go over next.
\begin{itemize}
\item The first one is given by the \textbf{compatibility condition}
\begin{align}
\label{eqinitial001}
\hspace*{-10pt} \sigma_0(\xi) \!=\! \frac{\mu}{2} \, (B_0 - B_0^{-1})
+ \frac{\lambda}{2} \, \log \det (B_0) \cdot \id, \quad
B_0(\xi) \! := \! \DD \varphi_0(\varphi_0^{-1}(\xi;t)) \,
\DD \varphi_0^T(\varphi_0^{-1}(\xi;t)),
\end{align}
that prescribes the initial values for $\sigma(\xi,0)$.
Here we need to assume $\varphi_0(x) \in H^2(\Omega)$ implying \newline
$\sigma_0,\, B_0\in W^{1,2}(\Omega_{\xi})$, which is a necessary condition for
the application of the Schauder fixed point theorem.
\item The second one is another \textbf{compatibility condition} given by the assumption of an %\deleted{\red{weak?}} 
 ``initial equilibrium''
\begin{align}
\Div_{\xi} \, \sigma_0(\xi) = f(\xi,0) \qquad \iff \qquad \Div_x S_1(\DD \varphi_0(x)) = \widetilde f(x,0).
\end{align}
\item Moreover, as stated in $\eqref{M1}_4$ the boundary values of the diffeomorphism $\varphi(x,t)$ are assumed to be constant in time,
which can be seen from
\begin{align}
v(\xi,t) := \ddt \varphi(\varphi^{-1}(\xi;t),t) = \ddt \varphi(x,t) = 0
\qquad \text{on} \quad \Sigma_T \quad \iff \quad \varphi(x,t) \equiv
\varphi(x) \quad \forall \, t \in [0,T).
\end{align}
Hence on all of $\Sigma_T$ we have
$\varphi(x,t)\equiv\varphi(x,0)=:\varphi_0(x)$ which is assumed to be a
known quantity.
\item Lastly, it will sometimes be convenient to assume $v_0(\xi)\equiv0$ on
$\Omega_\xi$, meaning that the system is supposed to be in an equilibrium state at
$t=0$. 
\end{itemize}
\subsubsection{Well-posedness of a sub-problem}
We end this section by proving the well-posedness of the system \eqref{M1} in
terms of the velocity $v$ \textbf{at a given history} of $\sigma, W$ and $f$.
For fixed $t$ and given $\sigma=\sigma(\xi,t) \in H^1(\Omega_\xi;\,\SYM)$,
$W(\xi, t)\in H^1(\Omega_{\xi}; \, \mathfrak{so}(3))$, the left hand side of
$\eqref{M1}_2$ defines a linear operator
$\mathbb{A}(\sigma):H_0^1(\Omega_\xi;\,\R^3) \to H^{-1}(\Omega_\xi;\,\R^3)$ by
\begin{equation}
\label{Adef}
\mathbb{A}(\sigma).v:=-\Div_{\xi} \left[\H^{\ZJ}(\sigma).D\right].
\end{equation}
We wish to point out that, in the present setting, this result can be alternatively approached by the Lax-Milgram theorem. However, for the purposes of our follow-up \cite{blesgen2024} and to be ready to face more general nonlinear hyperelasticity, it is useful to directly embed the set-up into that of monotone operators.
\begin{Theorem}
\label{lem1}
Let $t\in(0,T)$ and $\sigma(\cdot,t)\in H^1(\Omega_\xi;\,\SYM)$,
$W(\cdot,t)\in H^1(\Omega_{\xi}; \mathfrak{so}(3))$ be given and assume
$\H^{\ZJ}(\sigma)$ satisfies positive definiteness (\ref{H1}). Then the
equation~(\ref{Adef}) defines a strictly monotone operator
$v\mapsto \mathbb{A}(\sigma).v$ on $H_0^1(\Omega_\xi;\,\R^3)$. For
$g\in L^2(\Omega_T)$ and a.e.~$t\in(0,T)$, there exists a unique solution
$v\in H_0^1(\Omega_\xi;\,\R^3)$ to the vector-valued equation
\begin{equation}
\label{M3fix}
- \Div_\xi \left[\H^{\ZJ}(\sigma(\xi,t)).D\right]
\;=\; g(\xi,t) \qquad \mbox{\text{in} } \quad \Omega_\xi,
\end{equation}
where $g$ is defined in \eqref{eqdefofg} and assumed to be given.
\end{Theorem}
%%We recall that a second-order operator
%%$Lu(\xi):=-\sum_{i,j=1}^3\partial_i(a_{ij}(\xi)\paj u)$ is
%%{\it uniformly elliptic} if there exist constants $0<c_1\le C_1$ with
%%\[ c_1|w|^2\le\sum_{i,j=1}^3a_{ij}(\xi)w_iw_j\le C_1|w|^2\qquad
%%\mbox{for all }w\in\R^3,\,\xi\in\Omega_\xi. \]
%%If $a_{ij}$ is continuous and $\overline{\Omega_\xi}$ is compact, the
%%upper bound is naturally satisfied.
%%In comparison, the condition (\ref{H1}) provides only the lower bound for 
%%the uniform ellipticity of $A(\sigma)$ on $\Hsym$.
\begin{proof}
In order to show that (\ref{Adef}) defines a monotone operator
$v \mapsto \mathbb{A}(\sigma).v$, we need to verify two properties.

The first observation is that $\mathbb{A}(\sigma)$ is {\bf hemi-continuous},
i.e.~for arbitrary $v_1,\,v_2\in H_0^1(\Omega_\xi;\,\R^3)$ and any \break
$w\in H_0^1(\Omega_\xi;\,\R^3)$, it holds
\begin{align}
s \mapsto \big\llangle w, \mathbb{A}(\sigma).((1-s) \, v_1 + s \, v_2)
\big\rrangle \in C^0([0,1];\,\R).
\end{align}
Here and below, $\llangle w,w^*\rrangle$ denotes the duality pairing for
$w\!\in\! H_0^1(\Omega_\xi;\,\R^3)$, $w^*\!\in\! H^{-1}(\Omega_\xi;\,\R^3)$. 

Secondly, $v\mapsto \, \mathbb{A}(\sigma).v$ is {\bf monotone}.
For $v_1,v_2\in H_0^1(\Omega_\xi;\,\R^3)$, after integrating by parts
and using the positive definiteness of $\H^{\ZJ}(\sigma)$ from (\ref{H1}),
we have
\begin{equation}
\begin{alignedat}{2}
\big\llangle v_1-v_2,&\,\mathbb{A}(\sigma).v_1-\mathbb{A}(\sigma).v_2\big\rrangle
= -\io \langle (v_1-v_2) , \Div_\xi \! \left[(\H^{\ZJ}(\sigma))
.(D_1 - D_2)\right] \rangle_{\R^3} \dxi\\
&= \int_{\Omega_\xi} \langle \H^{\ZJ}(\sigma).(D_1 - D_2), \DD v_1 - \DD v_2 \rangle \dif \xi
= \io \big\langle \H^{\ZJ}(\sigma). (D_1 - D_2), D_1 - D_2\big\rangle\dxi\\
&\ge c_0\io\|D_1 - D_2\|^2\dxi  = c_0 \io \norm{\sym \, (\DD_\xi v_1 - \DD_\xi v_2)}^2 \dif \xi
% \;=\; c_0\io|\DD(v_1-v_2)|^2\dxi.
\;\ge\; \frac{c_0}{2}\|\DD_\xi(v_1-v_2)\|_{L^2(\Omega_\xi)}^2,
\end{alignedat}
\end{equation}
where Korn's inequality is used in the last line
(cf.~Gmeineder et al.~\cite{DieningGmeineder2024,gmeineder2023, gmeineder2024},
Lewintan et al.~\cite{Lewintan2021, Lewintan2021b, Lewintan2021c} and
Neff et al.~\cite{Neff2015}) and we write $D_i=\sym\, (\DD_\xi v_i)$, $i=1,2$.
By the Poincar{\'e} inequality, we hence find for $v_1\not=v_2$
\begin{align}
\big\llangle v_1-v_2,\mathbb{A}(\sigma).v_1-\mathbb{A}(\sigma).v_2\big\rrangle>0 \, ,
\end{align}
which is the strict monotonicity of $v \mapsto \mathbb{A}(\sigma).v$. 

From the hemi-continuity and the monotonicity of $\mathbb{A}(\sigma)$,
the existence of a solution to (\ref{M3fix}) follows from the theory of
monotone operators, see, e.g.~\cite{Zeidler2b}. The uniqueness of the solution
is a direct consequence of the strict monotonicity of $\mathbb{A}(\sigma)$.
\end{proof}
\begin{Remark}
\label{rem6}
Theorem~\ref{lem1} and its proof reveal that the major symmetry of $\H^{\ZJ}(\sigma)$
is not required and that instead the minor symmetry can be the right
framework for solving the divergence equation \eqref{M1}.
\end{Remark}
\begin{Remark}
\label{rem7}
Advanced elliptic regularity for PDE-systems in divergence form without major
symmetry can be found in Haller-Dintelmann et al.~\cite{Haller2019}.
If, in addition, we know that $\H^{\ZJ}(\sigma)$ has major symmetry, as is the case for 
equation $\eqref{M1}_2$ for given $\sigma(\xi,t), W(\xi,t), f(\xi,t)$ has variational
structure. Indeed, the corresponding minimization problem is given by
\begin{align}
\io \langle \H^{\ZJ}(\sigma(t)) .  \underbrace{\sym \, \DD_{\xi} v}_{=:D},
\sym \, \DD_{\xi} v \rangle - \langle g(\xi, t) ,  v(\xi,t) \rangle \dif \xi
\quad \longrightarrow \quad \min.v, \qquad v \in H_0^1(\Omega_{\xi};\R^3).
\end{align}
Again, existence and uniqueness follow from Korn's inequality.
The advantage of using this framework is that we may use full elliptic
regularity for the velocity field $v$ if wanted.
\end{Remark}
\subsection{Conditional determination of the diffeomorphism $\varphi(x,t)$ from the \newline solution $(\sigma,v)$ for $0<t<T^*$}
Assume that the spatial velocity $v(\xi,t)$ is regular enough to define
$\varphi(x,t)$ as the (unique) solution of the characteristic system
\begin{align}
\label{eqgewdgl01}
\partial_t \varphi(x,t) = v(\varphi(x,t),t), \qquad \quad \varphi(x,0) = \varphi_0(x).
\end{align}
This is for instance the case if $v(\xi,t)$ is Lipschitz continuous in $\xi$, which, for a
domain $\Omega_\xi$ that is smooth enough (Lipschitz boundary suffices), is
equivalent to $v(\cdot,t) \in W^{1,\infty}(\Omega_{\xi})$ and continuous
in $t$.
%\blue{(der Existenzsatz liefert $L^2$ in der Zeit und $W^{1,2}$ im Raum) Vor Einreichen entfernen $\dots$}. 
Then the Picard-Lindelöf theorem can be used to prove unique
solvability of \eqref{eqgewdgl01} and we may use the reconstructed deformation
$\varphi(x,t)$ to determine the quantities $F = \DD \varphi(x,t)$ and
$B = F \, F^T$. 

\medskip
\noindent Reconstruction of the constitutive law is possible by using the initial condition
\begin{align}
\sigma_0(\xi) = \frac{\mu}{2} \, (B_0 - B_0^{-1})
+ \frac{\lambda}{2} \, \log \det (B_0) \cdot \id, \qquad B_0(\xi) :=
\DD \varphi_0(\varphi_0^{-1}(\xi;t)) \,
\DD \varphi_0^T(\varphi_0^{-1}(\xi;t)) \, .
\end{align}

\medskip
Furthermore, for a sufficiently smooth solution $\varphi(x,t)$ of
\eqref{eqgewdgl01} we then have with $\dot{J} = \tr(L) \, J = \tr(D) \, J$
and $J(x,t) = \det F(x,t)$:
\begin{align}
\ddt[\log J] = \frac{\dot{J}}{J} = \tr(D)
\end{align}
so that $\log J(t) = \log J(0) + \int_0^t \tr(D(s)) \dif s$ and
\begin{align}
J(t) = J(0) \cdot \exp\Big(\int_0^t \tr(D(s)) \dif s\Big),
\end{align}
showing that automatically
$
\det \DD \varphi(x,t) = J(x,t) > 0.
$
Therefore, local 
%\deleted{\red{(local weak)}} 
smooth solutions of the hypoelastic
problem satisfy automatically the local invertibility constraint.
\begin{figure}[h!]
\begin{center}
\begin{minipage}[h!]{0.9\linewidth}
\if\Bilder y
\centering
\includegraphics[scale=0.6]{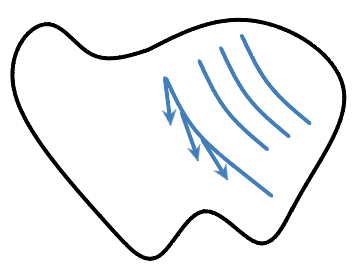}
% vorher (-35,115)
\put(-55,28){$\gamma(s)$}
\fi
\caption{Reconstructing the diffeomorphism $\varphi(x,t)$ by solving the characteristic system \newline $\partial_t \varphi(x,t) = v(\varphi(x,t),t)$ with $\varphi(x,0) = \varphi_0(x)$.}
\label{fig12}
\end{minipage}
\end{center}
\end{figure}
%
%\clearpage
%
% \blue{
% \noindent Bemerkung (conditional equilibrium): \\
% \\
% Vielleicht geht es konzeptionell nur so:
% 	\begin{itemize}
% 	\item Erst $v \in W^{1, \infty} \; \implies \; \varphi$ bestimmen ($\varphi$ sollte dabei lokal (also um den Anfangswert herum) $C^1$ sein).
% 	\item Mit $\varphi(x,t)$ $F$ und $B$ bestimmen sowie $W, L, D$.
% 	\item Compatible initial conditions $\varphi_0 \rightarrow \sigma_0$
% 		\begin{align}
% 		\partial_t \sigma + \sigma \, \underline W - \underline W \, \sigma = \underline{\H^{\ZJ}(\sigma)}. \underline D - \DD_\xi \sigma. \underline v.
% 		\end{align} 
% 	\item Zeige $\sigma = \sigma(B) = B - B^{-1} + \log \det B \cdot \id$ ist erfüllt.
% 	\item Jetzt mit Piola zurück zur Ref-Config, dort $\Div_x \frac{\dif}{\dif t} S_1 = \frac{\dif}{\dif t} \widetilde f$.
% 	\item Von hier jetzt ``integrieren'' erlaubt $\Div_x S_1 = \widetilde f \; \overset{\text{Piola}}{\iff} \; \Div_\xi \sigma = f$.
% 	\end{itemize}
% Aktuell scheint dieser Ansatz notwendig, um den Diffeomorphismus $\varphi(x,t)$ zu ermitteln und darüber die Gleichung $\Div_\xi \sigma(\xi,t) = f(\xi,t)$ als gelöst anzusehen (was ungünstig ist, da der Existenzsatz weniger liefert als das hier).
% }

\section{Conclusion}

Any hyperelastic isotropic nonlinear elasticity formulation can be written as a spatial rate-form equilibrium problem using objective stress rates and involving exclusively the Cauchy stress $\sigma$ and the spatial velocity $v$. For this, the Piola transformation is needed. Here, for simplicity, we considered a Cauchy-elastic formulation (not necessarily hyperelastic) for which all appearing tensor quantities can be made explicit and for which their properties are directly determined. In addition, we have shown that there emerges a subproblem of elliptic type that admits a unique solution based on Korn's inequality.\\
In the follow-up paper \cite{blesgen2024}
%part II of this work 
we will use the presented exceptional structure to provide a local existence result by applying a parabolic regularisation, using the properties of a differential
Sylvester equation together with the positive definiteness of
$\H^{\ZJ}(\sigma)$ and passing to the limit together with
a Schauder fixed-point argument. 

\medskip
\begin{comment}
\red{
\begin{align}
\frac{\DD^{\ZJ}}{\DD t}[\sigma](\xi,t) &= \H^{\ZJ}(\sigma(\xi,t)). D, \qquad
D = \sym \, \DD_{\xi} v && \qquad \text{in} \quad \Omega_T, \notag \\
-\div_\xi \left[\H^{\ZJ}(\sigma(B(\varphi^{-1}(\xi,t))),t). D \right] &=
-\ddt f(\xi,t) + \div[W \, \sigma - \sigma \, W](\xi,t) && \qquad \text{in}
\quad \Omega_T, \notag \\
v(\xi,0) &= v_0(\xi), \qquad \sigma(\xi,0) = \sigma_0(\xi) &&\qquad \text{in}
\quad \Omega_\xi, \notag \\
v(\xi,t) &= 0 &&\qquad \text{on} \quad \Sigma_T, \notag \\
\label{lasteq001}
\sigma(B) &= \frac{\mu}{2} \, (B- B^{-1}) + \frac{\lambda}{2} \,
\log \det B \cdot \id, \\
\partial_t \DD \varphi(x,t) &= \DD_{\xi}v(\varphi(x,t)) \cdot \DD \varphi(x,t),
\notag \\
\DD \varphi(x,0) &= \DD \varphi_0(x), \qquad \varphi(x,0) = \varphi_0(x). \notag
\end{align}
Blesgen: Eventuell müssen wir alles zusammen betrachten, um die conditional
positive definiteness von $\H^{\ZJ}(\sigma)$ aufzulösen? Dazu muss man aber
immer $B$ kennen, daher muss $\eqref{lasteq001}_6, \; \eqref{lasteq001}_7$
dazu.}
\end{comment}
%\\
%
\bigskip\par
\begingroup
\footnotesize
\noindent \textbf{Acknowledgement:} Patrizio Neff is grateful for discussions with
Davide Bigoni (University of Trento, Italy) and Robin J. Knops
(University of Edinburgh, Scotland) regarding failure of local uniqueness and
loss of stability in nonlinear elasticity for non rank-one convex formulations
on the occasion of the Euromech Colloquium ``630 Nonlinear Elasticity:
Modelling of multi-physics and applications - a Euromech/ICMS colloquium
celebrating the 80th birthday of Prof. Ray Ogden FRS''
(Edinburgh, UK, 25 March 2024 - 28 March 2024). \\
Patrizio Neff also acknowledges discussions with S. N. Korobeynikov
(Lavrentyev Institute of Hydrodynamics of Russian Academy of Science, Novosibirsk)
on the correct formulation of the equilibrium equations in the current configuration.
\endgroup
\begingroup
\footnotesize

\bibliographystyle{plain} %plain
\bibliography{references}
\endgroup

\begin{appendix}
\section{Appendix}
\subsection{Notation} \label{appendixnotation}
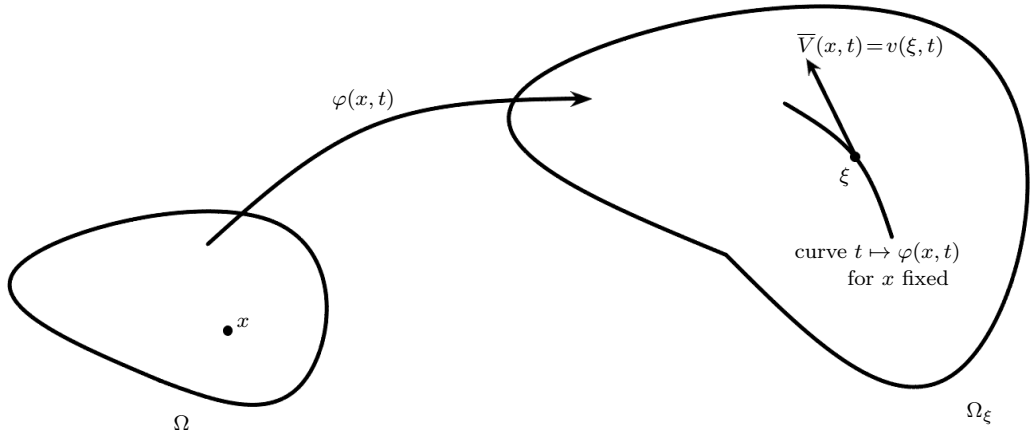
\begin{figure}[t]
		\begin{center}		
		\begin{minipage}[h!]{0.8\linewidth}
			\centering
			\hspace*{-40pt}
			\includegraphics[scale=0.4]{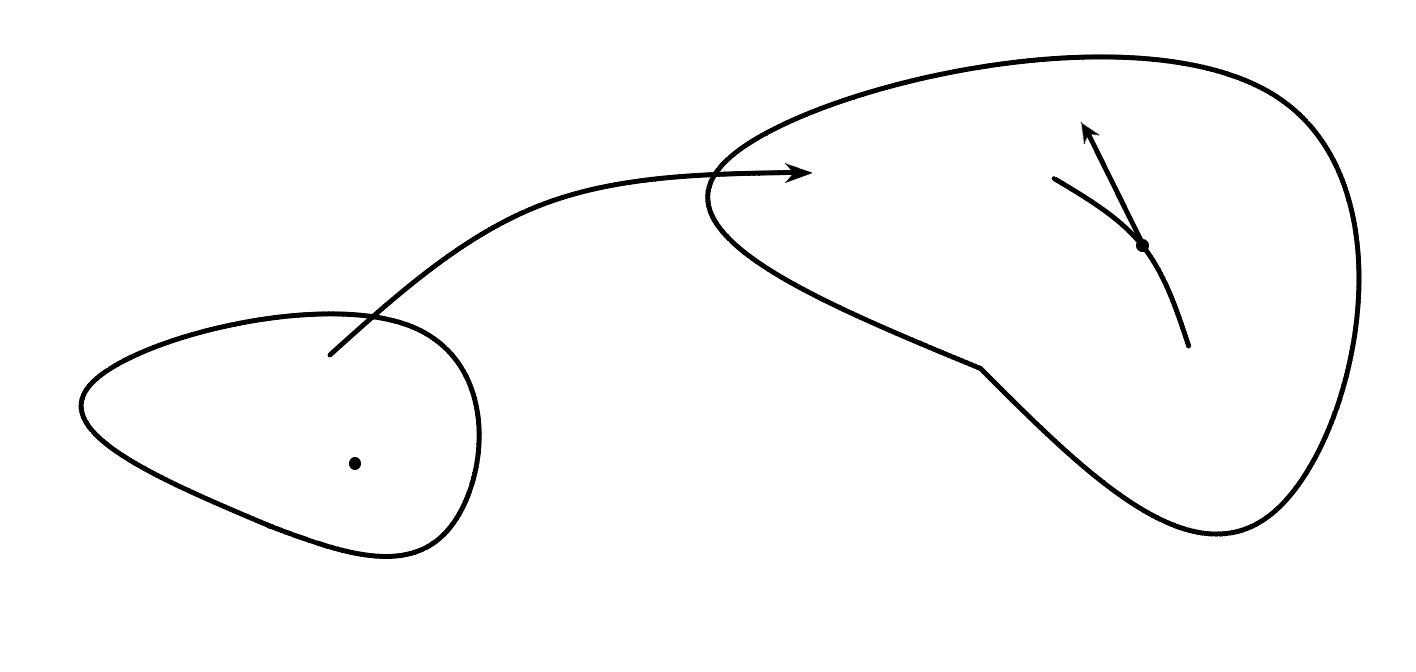}
			\put(-40,30){\footnotesize $\Omega_\xi$}
			\put(-340,25){\footnotesize $\Omega$}
			\put(-316,64){\footnotesize $x$}
			\put(-280,148){\footnotesize $\varphi(x,t)$}
			\put(-104,168){\footnotesize $\overline V(x,t) \!=\! v(\xi,t)$}
			\put(-88,119){\footnotesize $\xi$}
			\put(-105,90){\footnotesize curve $t \mapsto \varphi(x,t)$}
			\put(-85,80){\footnotesize  for $x$ fixed}
		\end{minipage} 
		\caption{Illustration of the deformation $\varphi(x,t): \Omega \to \Omega_{\xi}$ and the velocity $\overline V(x,t) = v(\xi,t)$.}
		\label{yfig1}
		\end{center}
	\end{figure}

\textbf{The deformation $\varphi(x,t)$, the material time derivative $\frac{\DD}{\DD t}$ and the partial time derivative $\partial_t$} \\
\\
In accordance with \cite{Marsden83} we agree on the following convention regarding an elastic deformation $\varphi$ and time derivatives of material quantities:

Given two sets $\Omega, \Omega_{\xi} \subset \R^3$ we denote by $\varphi: \Omega \to \Omega_{\xi}, x \mapsto \varphi(x) = \xi$ the deformation from the \emph{reference configuration} $\Omega$ to the \emph{current configuration} $\Omega_{\xi}$. A \emph{motion} of $\Omega$ is a time-dependent family of deformations, written $\xi = \varphi(x,t)$. The \emph{velocity} of the point $x \in \Omega$ is defined by $\overline{V}(x,t) = \partial_t \varphi(x,t)$ and describes a vector emanating from the point $\xi = \varphi(x,t)$ (see also Figure \ref{yfig1}). Similarly, the velocity viewed as a function of $\xi \in \Omega_{\xi}$ is denoted by $v(\xi,t)$. 

Considering an arbitrary material quantity $Q(x,t)$ on $\Omega$, equivalently represented by $q(\xi,t)$ on $\Omega_\xi$, we obtain by the chain rule for the time derivative of $Q(x,t)$
	\begin{align}
	\frac{\DD}{\DD t}q(\xi,t) \colonequals \frac{\dif}{\dif t}[Q(x,t)] = \DD_\xi q(\xi,t).v + \partial_t q(\xi,t) \, .
	\end{align}
Since it is always possible to view any material quantity $Q(x,t) = q(\xi,t)$ from two different angles, namely by holding $x$ or $\xi$ fixed, we agree to write
	\begin{itemize}
	\item $\dot q \colonequals \dd \frac{\DD}{\DD t}[q]$ for the material (substantial) derivative of $q$ with respect to $t$ holding $x$ fixed and
	\item $\partial_t q$ for the derivative of $q$ with respect to $t$ holding $\xi$ fixed.
	\end{itemize}
For example, we obtain the velocity gradient $L := \DD_\xi v(\xi,t)$ by
	\begin{align}
	L = \DD_\xi v(\xi,t) = \DD_\xi \overline V(x,t) \overset{\text{def}}&{=} \DD_\xi \frac{\dif}{\dif t} \varphi(x,t) = \DD_{\xi} \partial_t \varphi(\varphi^{-1}(\xi,t),t) = \partial_t \DD \varphi(\varphi^{-1}(\xi,t),t) \, \DD_\xi \big(\varphi^{-1}(\xi,t)\big) \notag \\
&=  \partial_t \DD \varphi(\varphi^{-1}(\xi,t),t) \, (\DD \varphi)^{-1}(\varphi^{-1}(\xi,t),t) = \dot F(x,t) \, F^{-1}(x,t) = L \, ,
	\end{align}
where we used that $\partial_t = \frac{\dif}{\dif t} = \frac{\DD}{\DD t}$ are all the same, if $x$ is fixed. \\
\\
As another example, when determining a corotational rate $\frac{\DD^{\circ}}{\DD t}$ we write
	\begin{align}
	\frac{\DD^{\circ}}{\DD t}[\sigma] = \frac{\DD}{\DD t}[\sigma] + \sigma \, \Omega^{\circ} - \Omega^{\circ} \, \sigma = \dot \sigma + \sigma \, \Omega^{\circ} - \Omega^{\circ} \, \sigma \, .
	\end{align}
However, if we solely work on the current configuration, i.e.~holding $\xi$ fixed, we write $\partial_t v$ for the time-derivative of the velocity (or any quantity in general). \\
\\
\textbf{Inner product} \\
\\
For $a,b\in\R^n$ we let $\langle {a},{b}\rangle_{\R^n}$  denote the scalar product on $\R^n$ with associated vector norm $\norm{a}_{\R^n}^2=\langle {a},{a}\rangle_{\R^n}$. We denote by $\R^{n\times n}$ the set of real $n\times n$ second-order tensors, written with capital letters. The standard Euclidean scalar product on $\R^{n\times n}$ is given by
$\langle {X},{Y}\rangle_{\R^{n\times n}}=\tr{(X Y^T)}$, where the superscript $^T$ is used to denote transposition. Thus the Frobenius tensor norm is $\norm{X}^2=\langle {X},{X}\rangle_{\R^{n\times n}}$, where we usually omit the subscript $\R^{n\times n}$ in writing the Frobenius tensor norm. The identity tensor on $\R^{n\times n}$ will be denoted by $\id$, so that $\tr{(X)}=\langle {X},{\id}\rangle$. \\
\\
\noindent \textbf{Frequently used spaces} 
	\begin{itemize}
	\item $\Sym(n), \rm \Sym^+(n)$ and $\Sym^{++}(n)$ denote the symmetric, positive semi-definite symmetric and positive definite symmetric second order tensors, respectively. Note that $\Sym^{++}(n)$ is considered herein only as an algebraic subset of $\Sym(n)$, not endowed with a Riemannian geometry \cite{fiala2009, fiala2016, fiala2020objective, kolev2024objective}.
	\item ${\rm GL}(n)\colonequals\{X\in\R^{n\times n}\;|\det{X}\neq 0\}$ denotes the general linear group.
	\item ${\rm GL}^+(n)\colonequals\{X\in\R^{n\times n}\;|\det{X}>0\}$ is the group of invertible matrices with positive determinant.
	\item ${\rm SL}(n)\colonequals\{X\in {\rm GL}(n)\;|\det{X}=1\}$.
	\item $\mathrm{O}(n)\colonequals\{X\in {\rm GL}(n)\;|\;X^TX=\id\}$.
	\item ${\rm SO}(n)\colonequals\{X\in {\rm GL}(n,\R)\;|\; X^T X=\id,\;\det{X}=1\}$.
	\item $\mathfrak{so}(3)\colonequals\{X\in\mathbb{R}^{3\times3}\;|\;X^T=-X\}$ is the Lie-algebra of skew symmetric tensors.
	\item $\mathfrak{sl}(3)\colonequals\{X\in\mathbb{R}^{3\times3}\;|\; \tr({X})=0\}$ is the Lie-algebra of traceless tensors.
	\item The set of positive real numbers is denoted by $\R_+\colonequals(0,\infty)$, while $\overline{\R}_+=\R_+\cup \{\infty\}$.
	\end{itemize}
\textbf{Frequently used tensors}
	\begin{itemize}
	\item $F = \DD \varphi(x,t)$ is the Fréchet derivative (Jacobian) of the deformation $\varphi(\,,t) : \Omega \to \Omega_{\xi} \subset \R^3$. $\varphi(x,t)$ is usually assumed to be a diffeomorphism at every time $t \ge 0$ so that the inverse mapping $\varphi^{-1}(\,,t) : \Omega_{\xi} \to \Omega$ exists.
	\item $C=F^T \, F$ is the right Cauchy-Green strain tensor.
	\item $B=F\, F^T$ is the left Cauchy-Green (or Finger) strain tensor.
	\item $U = \sqrt{F^T \, F} \in \Sym^{++}(3)$ is the right stretch tensor, i.e.~the unique element of ${\rm Sym}^{++}(3)$ with $U^2=C$.
	\item $V = \sqrt{F \, F^T} \in \Sym^{++}(3)$ is the left stretch tensor, i.e.~the unique element of ${\rm Sym}^{++}(3)$ with $V^2=B$.
	\item $\log V = \frac12 \, \log B$ is the spatial logarithmic strain tensor or Hencky strain.
        \item We write $V = Q$ diag($\lambda_1, \lambda_2, \lambda_3$) $Q^T$, where $\lambda_i \in \R_+$ are the principal stretches.
	\item $L = \dot F \, F^{-1} = \DD_\xi v(\xi)$ is the spatial velocity gradient.
	\item $v = \frac{\DD}{\DD t} \varphi(x, t)$ denotes the Eulerian velocity.
	\item $D = \sym \, L$ is the spatial rate of deformation, the Eulerian rate of deformation tensor (or stretching).
	\item $W = \sk \, L$ is the vorticity tensor or spin.
	\item We also have the polar decomposition $F = R \, U = V R \in {\rm GL}^+(3)$ with an orthogonal matrix $R \in \OO(3)$ (cf. Neff et al.~\cite{Neffpolardecomp}), see also \cite{LankeitNeffNakatsukasa,Neff_Nagatsukasa_logpolar13}.
	\end{itemize}
\noindent \textbf{Calculus with the material derivative} \\
\\
Consider the spatial Cauchy stress
	\begin{align}
	\label{eqmat01}
	\sigma(\xi,t) \colonequals \Sigma(B) = \Sigma(F(x,t) \, F^T(x,t)) = \Sigma(F(\varphi^{-1}(\xi,t),t) \, F^T(\varphi^{-1}(\xi,t),t)) \, .
	\end{align}
Then, on the one hand we have for the material derivative
	\begin{align}
	\label{eqmat02}
	\frac{\DD}{\DD t}[\sigma] = \DD_\xi \sigma(\xi,t).v(\xi,t) + \partial_t \sigma(\xi,t) \cdot 1
	\end{align}
and on the other hand equivalently
	\begin{align}
	\label{eqmat03}
	\frac{\DD}{\DD t}[\sigma] &= \frac{\DD}{\DD t}[\Sigma(F(x,t) \, F(x,t)^T)] \overset{(1)}{=} \frac{\dif}{\dif t}[\Sigma(F(x,t) \, F^T(x,t))] \\
	\overset{\substack{\text{standard} \\ \text{chain rule}}}&{=} \DD_B \Sigma(F(x,t) \, F^T(x,t)). \frac{\dif}{\dif t}[(F(x,t) \, F^T(x,t))] = \DD_B \Sigma(F(x,t) \, F^T(x,t)).(\dot F \, F^T + F \, \dot{F}^T) \nonumber\\
	&= \DD_B \Sigma(F(x,t) \, F^T(x,t)).[\dot F \, F^{-1} \, F \, F^T + F \, F^T \, F^{-T} \, \dot{F}^T] =\DD_B \Sigma(F(x,t) \, F^T(x,t)).[L \, B + B \, L^T] \, . \nonumber
	\end{align}
In $\eqref{eqmat03}_{1}$ we have used the fact that there is already a material representation which allows to set $\frac{\DD}{\DD t} = \frac{\dif}{\dif t}$. Of course, \eqref{eqmat02} is equivalent to \eqref{eqmat03}. From the context it should be clear which representation of $\sigma$ (referential or spatial) we are working with and by abuse of notation we do not distinguish between $\sigma$ and $\Sigma$. \\
\\
The same must be observed when calculating with corotational derivatives
	\begin{align}
	\label{eqmat04}
	\frac{\DD^{\circ}}{\DD t}[\sigma] = \frac{\DD}{\DD t}[\sigma] + \sigma \, \Omega^{\circ} - \Omega^{\circ} \, \sigma, \qquad \Omega^{\circ} = \frac{\DD}{\DD t}Q^{\circ}(x,t) \, (Q^{\circ})^T(x,t) = \frac{\dif}{\dif t} Q^{\circ}(x,t) \, (Q^{\circ})^T(x,t) \, .
	\end{align}
Here, we have
	\begin{align}
	\label{eqmat05}
	\frac{\DD^{\circ}}{\DD t}[\sigma] \overset{(\ast \ast)}&{=} Q^{\circ}(x,t) \, \frac{\DD}{\DD t}[(Q^{\circ})^T(x,t) \, \sigma \, Q^{\circ}(x,t)] \, (Q^{\circ})^T(x,t) \\
	&=Q^{\circ}(x,t) \, \left\{\frac{\DD}{\DD t}(Q^{\circ})^T(x,t) \, \sigma \, Q^{\circ}(x,t) + (Q^{\circ})^T(x,t) \, \frac{\DD}{\DD t}[\sigma] \, Q^{\circ}(x,t) + (Q^{\circ})^T(x,t) \, \sigma \, \frac{\DD}{\DD t}Q^{\circ}(x,t) \right\} \, (Q^{\circ})^T(x,t) \notag \\
	&=Q^{\circ}(x,t) \, \bigg\{\frac{\dif}{\dif t}(Q^{\circ})^T(x,t) \, \sigma \, Q^{\circ}(x,t) + (Q^{\circ})^T(x,t) \, \underbrace{\frac{\DD}{\DD t}[\sigma]}_{(\ast \ast \ast)} \, Q^{\circ}(x,t) + (Q^{\circ})^T(x,t) \, \sigma \, \frac{\dif}{\dif t}Q^{\circ}(x,t) \bigg\} \, (Q^{\circ})^T(x,t) \notag 
	\end{align}
and we can decide for $(\ast \ast \ast)$ to continue the calculus with \eqref{eqmat02} or \eqref{eqmat03}. In either case one has to decide viewing the functions as defined on the reference configuration $\Omega$ or in the spatial configuration $\Omega_\xi$.

In \eqref{eqmat05} we used $Q = Q(x,t)$ and $\Omega = \Omega(x,t)$. This means that the ``Lie-type'' representation $(\ast \, \ast)$ necessitates the definition of a reference configuration, so that we can switch between $\xi = \varphi(x,t)$ and $x$.

The interpretation $(\ast \, \ast)$ is most clearly represented for the Green-Naghdi rate $\frac{\DD^{\GN}}{\DD t}$, in which the spin \break $\Omega^{\GN} \colonequals \frac{\dif}{\dif t}R(x,t) \, R^T(x,t) = \dot R(x,t) \, R^T(x,t)$ is defined according to the polar decomposition $F = R \, U$ and in
	\begin{align}
	\frac{\DD^{\GN}}{\DD t}[\sigma] = \frac{\DD}{\DD t}[\sigma] + \sigma \, \Omega^{\GN} - \Omega^{\GN} \, \sigma = R \, \frac{\DD}{\DD t}[R^T \, \sigma \, R] \, R^T
	\end{align}
the term $[R^T \, \sigma \, R]$ is called \emph{corotational stress tensor} (cf.~\cite[p.~142]{Marsden83}). \\
\\
\noindent \textbf{Tensor domains} \\
\\
Denoting the reference configuration by $\Omega$ with tangential space $T_x \Omega$ and the current/spatial configuration by $\Omega_\xi$ with tangential space $T_\xi \Omega_\xi$ as well as $\varphi(x) = \xi$, we have the following relations (see also Figure \ref{yfig2}):

	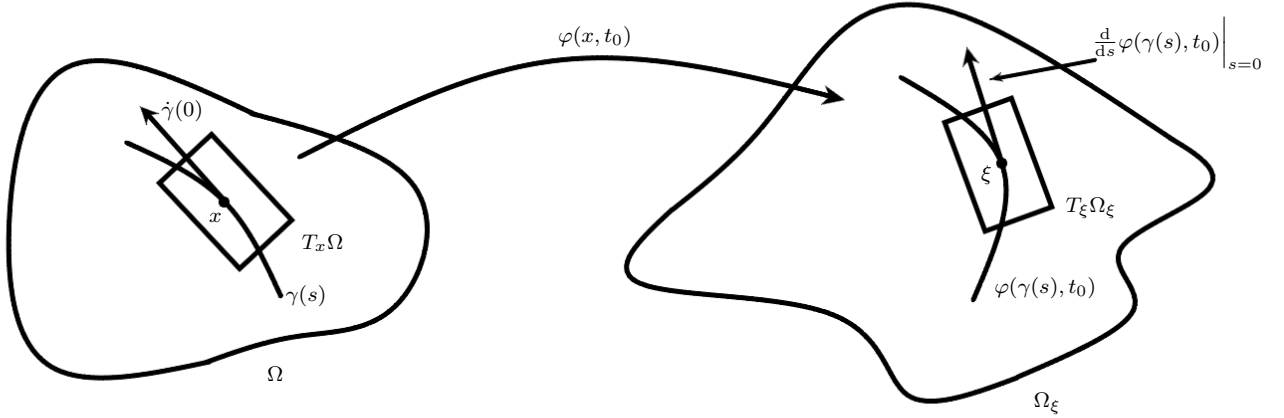
\begin{figure}[t]
		\begin{center}		
		\begin{minipage}[h!]{0.8\linewidth}
			\centering
			\hspace*{-80pt}
			\includegraphics[scale=0.5]{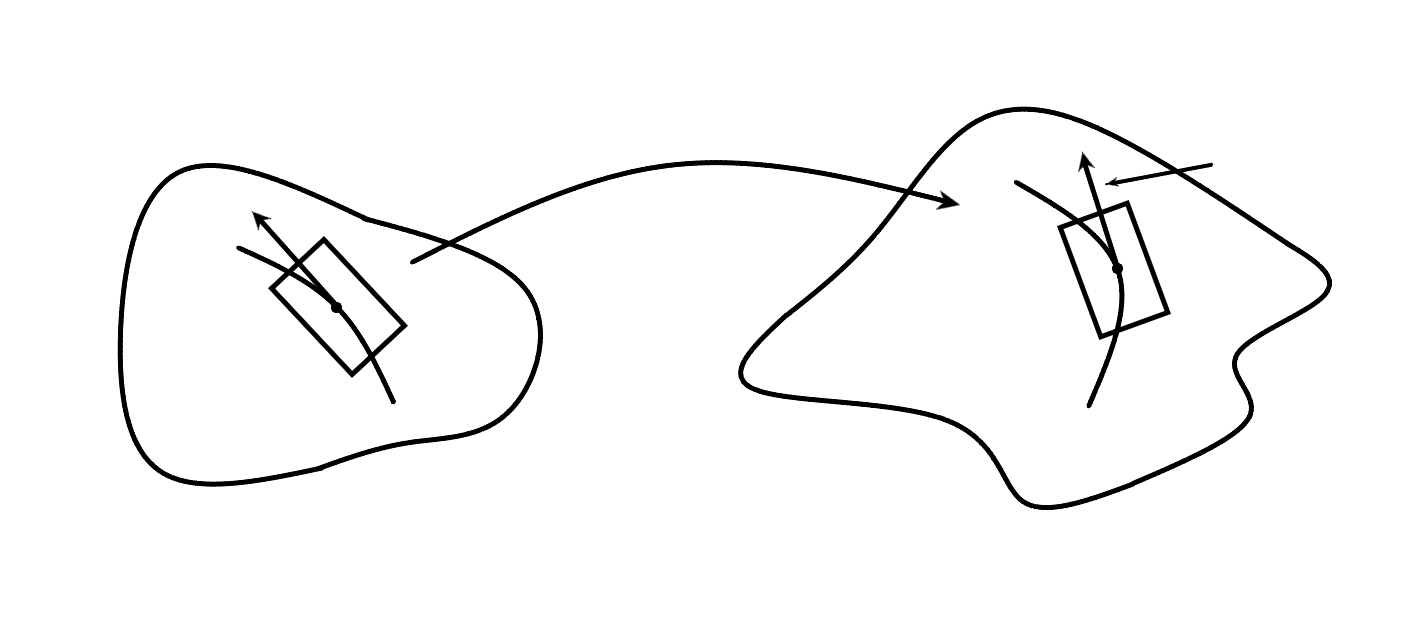}
			\put(-100,45){\footnotesize $\Omega_\xi$}
			\put(-390,55){\footnotesize $\Omega$}
			\put(-412,115){\footnotesize $x$}
			\put(-430,155){\footnotesize $\dot \gamma(0)$}
			\put(-377,105){\footnotesize $T_x \Omega$}
			\put(-383,85){\footnotesize $\gamma(s)$}
			\put(-280,183){\footnotesize $\varphi(x,t_0)$}
			\put(-120,131){\footnotesize $\xi$}
			\put(-78,181){\footnotesize $\frac{\dif}{\dif s}\varphi(\gamma(s),t_0)\bigg\vert_{s=0}$}
			\put(-88,119){\footnotesize $T_\xi \Omega_\xi$}
			\put(-115,88){\footnotesize $\varphi(\gamma(s),t_0)$}
		\end{minipage} 
		\caption{Illustration of the curve $s \mapsto \varphi(\gamma(s),t_0), \; \gamma(0) = x$ for a fixed time $t = t_0$ with vector field \break $s \mapsto \frac{\dif}{\dif s} \varphi(\gamma(s),t) \in T_\xi \Omega_\xi$.}
		\label{yfig2}
		\end{center}
	\end{figure}
	\begin{multicols}{3}
	\begin{itemize}
	\item $F\colon T_x \Omega \to T_\xi \Omega_\xi$
	\item $R\colon T_x \Omega \to T_\xi \Omega_\xi$
	\item $F^T\colon T_\xi \Omega_\xi \to T_x \Omega$
	\item $R^T\colon T_\xi \Omega_\xi \to T_x \Omega$
	\item $C = F^T \, F\colon T_x \Omega \to T_x \Omega$
	\item $B = F \, F^T\colon T_\xi \Omega_\xi \to T_\xi \Omega_\xi$
	\item $\sigma\colon T_\xi \Omega_\xi \to T_\xi \Omega_\xi$
	\item $\tau\colon T_\xi \Omega_\xi \to T_\xi \Omega_\xi$
	\item $S_2 = F^{-1} \, S_1 \colon T_x \Omega \to T_x \Omega$
	\item $S_1 = \DD_F \WW(F) \colon T_x \Omega \to T_\xi \Omega_\xi$
	\item $R^T \, \sigma \, R\colon T_x \Omega \to T_x \Omega$.
	\end{itemize}
	\end{multicols}

\noindent \textbf{The strain energy function $\WW(F)$} \\
\\
We are only concerned with rotationally symmetric functions $\WW(F)$ (objective and isotropic), i.e.
	\begin{equation*}
	\WW(F)={\WW}(Q_1^T\, F\, Q_2), \qquad \forall \, F \in {\rm GL}^+(3), \qquad  Q_1 , Q_2 \in {\rm SO}(3).
	\end{equation*}
\textbf{List of additional definitions and useful identities}
	\begin{itemize}
	\item For two metric spaces $X, Y$ and a linear map $L: X \to Y$ with argument $v \in X$ we write $L.v:=L(v)$. This applies to a second-order tensor $A$ and a vector $v$ as $A.v$ as well as to a fourth-order tensor $\C$ and a second-order tensor $H$ as $\C.H$.
	\item We denote the space of minor and major symmetric, positive definite fourth-order tensors $\C$ by $\Sym^{++}_4(6)$, i.e.~$\C \in \Sym^{++}_4(6)$ if and only if $\langle \C.D, D \rangle > 0$ for all $D \in \Sym(3) \! \setminus \! \{0\}$.
	\item By $\div(\cdot)$ we denote the common divergence of vectors, i.e.~$\div_\xi v=\sum_{i=1}^3\partial_{\xi_i}v_i$.
	\item By $\Div(\cdot)$ we denote the divergence of matrices, i.e.~for $A=(a_1,a_2,\ldots,a_n)\in\R^{n\times n}$, where $a_i\in \R^n$ are the rows of $A$, we have $\Div_\xi A=(\div_\xi\, a_1,\div_\xi\, a_2,\ldots,\div_\xi\, a_n)$.
	\item We define $J = \det{F}$ and denote by $\Cof(X) = (\det X)X^{-T}$ the cofactor of a matrix in ${\rm GL}^{+}(3)$.
	\item We define $\sym \, X = \frac12 \, (X + X^T)$ and $\sk X = \frac12 \, (X - X^T)$ as well as $\dev X = X - \frac13 \, \tr(X) \cdot \id$.
    %\blue{oder $\dev_3 X$}.
	\item For all vectors $\xi,\eta\in\R^3$ we have the tensor or dyadic product $(\xi\otimes\eta)_{ij}=\xi_i\,\eta_j$.
	\item $S_1=\DD_F \WW(F) = \sigma \cdot \Cof F$ is the non-symmetric first Piola-Kirchhoff stress tensor.
	\item $S_2=F^{-1}S_1=2\,\DD_C \widetilde{\WW}(C)$ is the symmetric second  Piola-Kirchhoff stress tensor.
	\item $\sigma=\frac{1}{J}\,  S_1\, F^T=\frac{1}{J}\,  F\,S_2\, F^T=\frac{2}{J}\DD_B \widetilde{\WW}(B)\, B=\frac{1}{J}\DD_V \widetilde{\WW}(V)\, V = \frac{1}{J} \, \DD_{\log V} \widehat \WW(\log V)$ is the symmetric Cauchy stress tensor.
	\item $\sigma = \frac{1}{J} \, F\, S_2 \, F^T = \frac{2}{J} \, F \, \DD_C \widetilde{\WW}(C) \, F^T$ is the ''\emph{Doyle-Ericksen formula}'' \cite{doyle1956}.
	\item For $\sigma: \Sym(3) \to \Sym(3)$ we denote by $\DD_B \sigma(B)$ with $\sigma(B+H) = \sigma(B) + \DD_B \sigma(B).H + o(H)$ the Fréchet-derivative. For $\sigma: \Sym^+(3) \subset \Sym(3) \to \Sym(3)$ the same applies. Similarly, for $\WW : \R^{3 \times 3} \to \R$ we have $\WW(X + H) = \WW(X) + \langle \DD_X \WW(X), H \rangle + o(H)$.
	\item $\tau = J \, \sigma = 2\, \DD_B \widetilde{\WW}(B)\, B $ is the symmetric Kirchhoff stress tensor.
	\item $\tau = \DD_{\log V} \widehat{\WW}(\log V)$ is the ``\emph{Richter-formula}'' \cite{richter1948isotrope, richter1949hauptaufsatze}.
	\item $\sigma_i =\dd\frac{1}{\lambda_1\lambda_2\lambda_3}\dd\lambda_i\frac{\partial g(\lambda_1,\lambda_2,\lambda_3)}{\partial \lambda_i}=\dd\frac{1}{\lambda_j\lambda_k}\dd\frac{\partial g(\lambda_1,\lambda_2,\lambda_3)}{\partial \lambda_i}, \ \ i\neq j\neq k \neq i$ are the principal Cauchy stresses (the eigenvalues of the Cauchy stress tensor $\sigma$), where $g:\mathbb{R}_+^3\to \mathbb{R}$ is the unique function of the singular values of $U$ (the principal stretches) such that $\WW(F)=\widetilde{\WW}(U)=g(\lambda_1,\lambda_2,\lambda_3)$.
	\item $\sigma_i =\dd\frac{1}{\lambda_1\lambda_2\lambda_3}\frac{\partial \widehat{g}(\log \lambda_1,\log \lambda_2,\log \lambda_3)}{\partial \log \lambda_i}$, where $\widehat{g}:\mathbb{R}^3\to \mathbb{R}$ is the unique function such that \\ \hspace*{0.3cm} $\widehat{g}(\log \lambda_1,\log \lambda_2,\log \lambda_3):=g(\lambda_1,\lambda_2,\lambda_3)$.
	\item $\tau_i =J\, \sigma_i=\dd\lambda_i\frac{\partial g(\lambda_1,\lambda_2,\lambda_3)}{\partial \lambda_i}=\frac{\partial \widehat{g}(\log \lambda_1,\log \lambda_2,\log \lambda_3)}{\partial \log \lambda_i}$ \, .
 %    \red{
 %    \item $T_{\Biot} = \DD_U \widetilde{\WW}(U)$ is the symmetric Biot stress tensor
	% % (symmetric only in the case of isotropy)
 %    \item $\sigma = \frac{1}{\det V} \, V (R \, T_{\Biot}\, R^T) $
 %    \item $T_{\Biot}^{i} = \frac{\partial g(\lambda_1,\lambda_2,\lambda_3)}{\partial \lambda_i}$ (in case of hyperelasticity)
 %    \item $\sigma_i = \frac{1}{\lambda_j \, \lambda_k} T_{\Biot}^{i} \; \left(= \; \lambda_i T_{\Biot}^i \; \text{for incompressibility} \right) $
 %    }
	\end{itemize}
\vspace*{2em}
\noindent \textbf{Conventions for fourth-order symmetric operators, minor and major symmetry}

\medskip\par\noindent%
For a fourth-order linear mapping $\C\colon\R^{3\times3} \to \R^{3\times3}$, we agree on the following convention.

\medskip\par\noindent%
We say that $\C$ has \emph{minor symmetry} if
\begin{align}
	\C.S \in \Sym(3) \qquad \forall \, S \in \Sym(3).
\end{align}
This can also be written in index notation as $\C_{ijkm} = \C_{jikm} = \C_{ijmk}$. 
%In general, $\C$ can be transformed into having minor symmetry by considering the mapping $X \to \sym(\C. \sym X)$, and $\C\colon \R^{3 \times 3} \to \R^{3 \times 3}$ is minor symmetric if and only if $\C.X = \sym(\C.\sym X)$. \textcolor{blue}{TODO: this is wrong ($\C=\operatorname{id}$)}

\medskip\par\noindent%
We say that $\C\colon \Sym(3) \to \Sym(3)$ has \emph{major symmetry} (or is \emph{self-adjoint}) if
\begin{align}
	\langle \C. T, S \rangle = \langle \C. S, T \rangle \qquad \forall \, T, S \in \Sym(3)\,.
\end{align}
Major symmetry in index notation is understood as $C_{ijkm} = C_{kmij}$.
For $\C\colon\Sym(3)\to\Sym(3)$, we define the adjoint operator $\C^T\colon\Sym(3)\to\Sym(3)$ via $\langle \C^T . T , S \rangle = \langle \C . S , T \rangle \ \forall S, T \in \Sym(3)$. Then $\sym \, \C \colonequals \frac{1}{2} \, (\C^T + \C)$ denotes the major-symmetric part of $\C$. The set of major (and minor) symmetric fourth-order tensor mappings $\Sym(3) \to \Sym(3)$ which are \emph{positive-definite} and \emph{positive-semidefinite} are denoted by $\Sym^{++}_4(6)$ and $\Sym^{+}_4(6)$, respectively, i.e.
\begin{alignat}{2}
	\label{eqposdef1}
	&\C \in \Sym^{++}_4(6)\qquad &&\iff \qquad \langle \C.H, H \rangle > 0 \qquad \forall \, H \in \Sym(3)
	\,,
	\\
	\label{eqpossemidef1}
	&\C \in \Sym^+_4(6) \qquad &&\iff \qquad \langle \C.H, H \rangle \ge 0 \qquad \forall \, H \in \Sym(3)
	\,.
\end{alignat}
By identifying $\Sym(3) \cong \R^6$, we can interpret $\C\colon\Sym(3) \to \Sym(3)$ as a linear mapping $\widetilde \C\colon \R^6 \to \R^6$. More specifically,
for $H \in \Sym(3) \cong \R^6$
%has the entries $H_{ij}$,
we can write, using the Mandel-notation,
\begin{align}
	\label{eqvec1}
	h \colonequals \text{vec}(H) = (H_{11}, H_{22}, H_{33}, \sqrt{2} \, H_{12}, \sqrt{2} \, H_{23}, \sqrt{2} \, H_{31}) \in \R^6
\end{align}
so that with $\widetilde{\C}. \text{vec}(H)\colonequals \text{vec}(\C.H)$ we have $\langle \C.H, H \rangle_{\Sym(3)} = \langle \widetilde \C.h, h \rangle_{\R^6}$.
Then for any fourth-order tensor $\C\colon \Sym(3) \to \Sym(3)$, we can also define $\textbf{sym} \, \C$ by $(\textbf{sym} \, \C).H = \text{vec}^{-1}((\sym \, \widetilde{\C}) . \text{vec}(H))$, implying
\begin{align}
	\langle \C.H, H \rangle_{\Sym(3)} = \langle \widetilde \C.h, h \rangle_{\R^6} = \langle (\sym  \, \widetilde \C). h, h \rangle_{\R^6} = \langle \textbf{sym} \, \C.H, H \rangle_{\Sym(3)} \qquad \forall \, H \in \Sym(3).
\end{align}
Major symmetry in these terms can be expressed as $\widetilde \C \in \Sym(6)$. However, we omit the tilde-operation and ${\bf sym}$ and write in short $\sym \, \C\in \Sym_4(6)$ if no confusion can arise. In the same manner we speak about $\det \C$ meaning $\det \widetilde \C$.

\subsection{Derivation of the induced fourth-order tangent stiffness tensor $\H^{\ZJ}(\sigma)$}
\label{sec:rate-formulations}
For the derivation of the induced
tangent stiffness tensor $\H^{\ZJ}(\sigma)$ for the Cauchy stress-strain law
\begin{align}
\sigma(B) = \frac{\mu}{2} \, (B - B^{-1}) + \frac{\lambda}{2} \, \log \det B
\cdot \id
\end{align}
we begin by splitting $\sigma(B)$ into three parts, namely
\begin{align}
\label{eqsigma2}
\sigma = \frac{\mu}{2} (B - B^{-1}) + \frac{\lambda}{2} \, \tr(\log B) \, \id
= \frac12 \left\{\underbrace{\mu \, (B-\id)}_{=: \sigma_1} +
\underbrace{\mu \, (\id-B^{-1})}_{=: \sigma_2}
+ \underbrace{\lambda \, \tr(\log B) \, \id}_{=:\sigma_3}\right\}.
\end{align}
Next, we calculate the Zaremba-Jaumann derivative, given by
\begin{align}
\frac{\DD^{\ZJ}}{\DD t}[\sigma] &= \frac{\DD}{\DD t}[\sigma] + \sigma \, W - W \, \sigma
\end{align}
for each of the three components. To do so, we use the easily verifiable identity $\frac{\DD}{\DD t}[B] = L \, B + B \, L^T$, yielding for $\sigma_1$
\begin{equation}
\label{eqsigma1}
\begin{alignedat}{2}
\frac{\DD}{\DD t} [\sigma_1] &= \mu \, \frac{\DD}{\DD t} [B- \id] = \mu \,
(L \, B + B \, L^T) = \mu \, (L \, (B-\id) + (B-\id) \, L^T + L + L^T) \\
&= (L \, \sigma_1 + \sigma_1 \, L^T) + 2 \, \mu \, D = (D+W) \, \sigma_1 + \sigma_1 \,
(D+W)^T + 2 \, \mu \, D \\
&= D \, \sigma_1 + \sigma_1 \, D + 2 \, \mu \, D + W\, \sigma_1 - \sigma_1 \, W \\
\iff \frac{\DD^{\ZJ}}{\DD t}[\sigma_1] &= \frac{\DD}{\DD t} [\sigma_1] + \sigma_1 \, W - W \,
\sigma_1 = 2 \, \mu \, D + (\sigma_1 \, D + D \, \sigma_1) \, .
\end{alignedat}
\end{equation}
Recalling the identities $\frac{\DD}{\DD t} [B^{-1}] = -B^{-1} \, \dot{B} \, B^{-1}$
and $B^{-1} = \id - \frac{1}{\mu} \, \sigma_2$, we obtain for $\sigma_2$
\begin{equation}
\begin{alignedat}{2}
\frac{\DD}{\DD t} [\sigma_2] &= \frac{\DD}{\DD t} [\mu \, (\id - B^{-1})] = -\mu \,
\frac{\DD}{\DD t} [B^{-1}] = \mu \, B^{-1} \, \dot{B} \, B^{-1} \\
&= \mu \, B^{-1} \, (L\, B + B \, L^T) \, B^{-1} = \mu \, (B^{-1} \, L + L^T \,
B^{-1}) \\
&= \mu (B^{-1} \, W + B^{-1} \, D -W \, B^{-1} + D \, B^{-1}) \\
&= \mu\left(\left(\id - \frac{1}{\mu} \, \sigma_2\right) \, D + D \,
\left(\id - \frac{1}{\mu} \, \sigma_2\right) + \left(\id - \frac{1}{\mu} \,
\sigma_2\right) \, W - W \left(\id - \frac{1}{\mu} \, \sigma_2\right)\right) \\
&= 2\mu \, D - (\sigma_2 \, D + D \, \sigma_2) + W \, \sigma_2 - \sigma_2 \, W\\
\implies \frac{\DD^{\ZJ}}{\DD t} [\sigma_2] &= 2\mu \, D -
(\sigma_2 \, D + D \, \sigma_2).
\end{alignedat}
\end{equation}
Also, since $\tr(\log B) \cdot \id = \log (\det B) \cdot \id$
(cf. Neff et al. \cite{Neff_Osterbrink_Martin_Hencky13, NeffGhibaLankeit}) and
$\sigma_3 \, W - W \, \sigma_3 = 0$, since
$\sigma_3 = g(t) \cdot \id, \, g(t) \in \R$, we have for $\sigma_3$
\begin{equation}
\begin{alignedat}{2}
\frac{\DD}{\DD t} [\log(\det B) \cdot \id] &= \frac{1}{\det B} \, \langle \Cof
B , \dot{B} \rangle \, \id = \frac{1}{\det B} \cdot \det B \cdot \langle B^{-1},
\dot{B} \rangle \cdot \id \\
&= \langle B^{-1}, L \, B + B \, L^T \rangle \cdot \id = ( \langle \id, L \rangle
+ \langle \id , L^T \rangle ) \cdot \id = 2 \, \tr(D) \cdot \id \\
\implies \frac{\DD^{\ZJ}}{\DD t}[\sigma_3] &= 2 \lambda \, \tr(D) \cdot \id.
\end{alignedat}
\end{equation}
Combining linearly these three expressions for the rates now yields for the
Cauchy stress $\sigma$
\begin{align}
\frac{\DD^{\ZJ}}{\DD t}[\sigma] &= \frac12 \,
\frac{\DD^{\ZJ}}{\DD t}[\sigma_1 + \sigma_2 + \sigma_3] = 2\mu \, D
+ \frac{1}{2}\{(\sigma_1 - \sigma_2) \, D + D \, (\sigma_1 - \sigma_2)\}
+ \lambda \, \tr(D) \cdot \id \notag \\
\label{eqallsigma}
&= 2\mu \, D + \frac{\mu}{2} \, \{(B + B^{-1}) \, D + D \, (B + B^{-1}) - 4 \,
D\} + \lambda \, \tr(D) \cdot \id \\
&= \frac{\mu}{2} \, \{B \, D + D \, B + B^{-1} \, D + D \, B^{-1}\}
+ \lambda \, \tr(D) \cdot \id =: \H^{\ZJ}(B).D = \H^{\ZJ}(\sigma).D, \notag
\end{align}
where the last equality holds since $\sigma: \Sym^{++}(3) \to \Sym(3), B \mapsto \sigma(B)$ is invertible.

\subsection{Salient properties of $\sigma(B) = \frac{\mu}{2} \, (B - B^{-1}) + \frac{\lambda}{2} \, (\log \det B) \cdot \mathbbm{1}$} 
\label{secA1}
The {\sl principal Cauchy stresses} $\sigma_i$ are the eigenvalues of $\sigma$. Since $\sigma$ is coaxial to $B$ (cf. Thiel et al \cite{thiel2019}) we have
	\begin{align}
	\sigma(Q^T \, B \, Q) = Q^T \, \sigma(B) \, Q = \diag(\sigma_1, \sigma_2, \sigma_3) = \sigma(\diag(\lambda_1^2, \lambda_2^2, \lambda_3^2)).
	\end{align}
Therefore, we may simply insert $B = \text{diag}(\lambda_1^2, \lambda_2^2, \lambda_3^2)$ into the constitutive law
	\begin{align}
	\label{eqappendix01}
	\sigma(B) = \frac{\mu}{2} \, (B - B^{-1}) + \frac{\lambda}{2} \, \log \det B \cdot \id
	\end{align}
and calculate the diagonal entries of $\sigma(B)$, yielding
	\begin{align}
	\sigma_i = \frac{\mu}{2} \, (\lambda_i^2 - \lambda_i^{-2}) + \lambda \, \log(\lambda_1 \, \lambda_2 \, \lambda_3),
	\end{align}
where $\lambda_i^2, \; i=1,2,3$ are the eigenvalues of $B$.
\subsubsection{Tension-extension inequality}
The principal Cauchy stresses $\sigma_i$ fulfill the tension-extension inequality $\frac{\partial \sigma_i}{\partial \lambda_i} > 0$. Indeed, for $\lambda_k, \lambda_j$ fixed, $k,j \neq i$,
	\begin{align}
	\frac{\partial \sigma_i}{\partial \lambda_i} = \mu \, \left(\lambda_i + \frac{1}{\lambda_i^3}\right) + \frac{\lambda}{\lambda_i} \overset{\lambda_i > 0}{>} 0, \qquad \mu, \lambda > 0.
	\end{align}
\subsubsection{Baker-Ericksen inequality}
The principal Cauchy stresses $\sigma_i$ fulfill the Baker-Ericksen inequality $(\sigma_i - \sigma_j) \, (\lambda_i - \lambda_j) > 0$: 
	\begin{equation}
	\begin{alignedat}{2}
	\sigma_i - \sigma_j = \frac{\mu}{2} \, \left(\lambda_i^2 - \lambda_j^2 - \left(\frac{1}{\lambda_i^2} - \frac{1}{\lambda_j^2}\right)\right) &= \frac{\mu}{2} \, \left\{\lambda_i + \lambda_j + \frac{\lambda_i + \lambda_j}{\lambda_i^2 \, \lambda_j^2}\right\} \, (\lambda_i - \lambda_j) \quad \text{for } i \neq j\\
	\implies \qquad (\sigma_i - \sigma_j) \, (\lambda_i - \lambda_j) &= \frac{\mu}{2} \, \left\{\lambda_i + \lambda_j + \frac{\lambda_i + \lambda_j}{\lambda_i^2 \, \lambda_j^2}\right\} \, (\lambda_i - \lambda_j)^2 > 0.
	\end{alignedat}
	\end{equation}
\subsubsection{Monotonicity of pressure: pressure-compression inequality}
The constitutive law \eqref{eqappendix01} fulfills the monotonicity of pressure, i.e. $\langle \sigma(\alpha \, \id) - \sigma(\beta \, \id), \alpha \, \id - \beta \, \id \rangle > 0$ for $\alpha, \beta > 0$ $\alpha \neq \beta$:
	\begin{align}
	\langle \sigma(\alpha \, \id) - \sigma(\beta \, \id), \alpha \, \id - \beta \, \id \rangle &= \Big\langle \left\{ \frac{\mu}{2} \, \left( \alpha - \frac{1}{\alpha} \right) \, \id + \frac{\lambda}{2} \cdot 3 \, \log \alpha \cdot \id \right\} - \left\{ \frac{\mu}{2} \, \left( \beta - \frac{1}{\beta} \right) \, \id + \frac{\lambda}{2} \cdot 3 \, \log \beta \cdot \id \right\}, \alpha \, \id - \beta \, \id \Big\rangle \notag \\
	&= \left[ \frac{\mu}{2} \cdot \left\{ \alpha - \beta + \frac{\alpha - \beta}{\alpha \, \beta} \right\} + \frac{3 \, \lambda}{2} \cdot (\log \alpha - \log \beta) \right] \cdot (\alpha - \beta) \, \langle \id , \id \rangle \\
	&= 3 \, \frac{\mu}{2} \cdot \left\{(\alpha - \beta)^2 + \frac{(\alpha - \beta)^2}{\alpha \, \beta} \right\} + 3 \, \frac{3 \, \lambda}{2} \cdot (\alpha - \beta) \cdot (\log \alpha - \log \beta) > 0. \notag
	\end{align}
The latter implies the pressure-compression inequality $\frac13 \, \tr(\sigma(\alpha \, \id)) \cdot (\alpha - 1) > 0$, by setting $\beta = 1$.
	\begin{figure}[h!]
		\begin{center}
		\begin{minipage}[h!]{0.9\linewidth}
			\if\Bilder y
			\centering
			\includegraphics[scale=0.45]{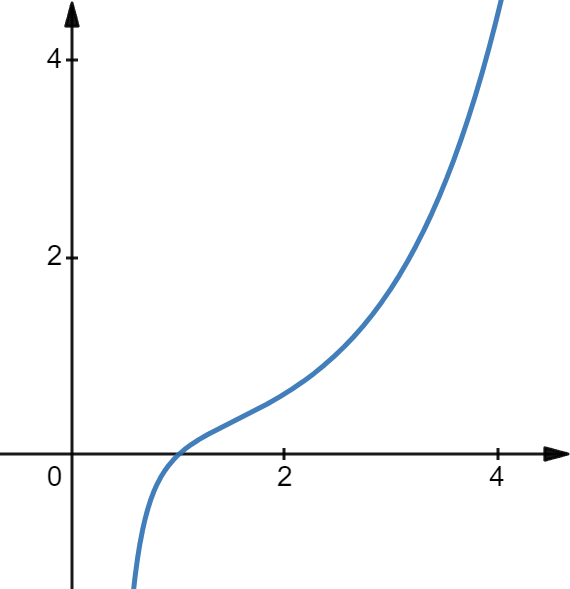}
			\put(-15,35){$\alpha$}
			\put(-135,35){$1$}
			\put(-160,180){$\frac{1}{3}\tr(\sigma(\alpha \,\id))$}
			\centering
			\fi
			\caption{Monotone volume-mean pressure $\alpha \mapsto \frac{1}{3}\tr(\sigma(\alpha \,\id))$. It is physically reasonable to expect that the mean pressure is an increasing function of the volumetric stretch for a stable material. }
			\label{graphic-dupa-k}
		\end{minipage}
	\end{center}
	\end{figure}
\subsubsection{Smooth invertibility of $B \mapsto \sigma(B) = \frac{\mu}{2} \,
(B - B^{-1}) + \frac{\lambda}{2} \, \log \det B \cdot \id$}
In an upcoming paper by Martin et al.~\cite{MartinVossGhibaNeff} it will be proven that
\begin{align}
\sigma: \Sym^{++}(3) \to \Sym(3), \qquad B \mapsto \sigma(B) = \frac{\mu}{2} \,
(B - B^{-1}) + \frac{\lambda}{2} \, \log \det B \cdot \id
\end{align}
is smoothly invertible. Hence, there exists an inverse function $\mathcal{F}^{-1}: \Sym(3) \to \Sym^{++}(3), \quad \sigma \mapsto \mathcal{F}^{-1}(\sigma)$.

\subsubsection{Positive definiteness}
The induced fourth-order tangent stiffness tensor $\H^{\ZJ}(\sigma)$ is positive definite for $\mu > 0, 3\, \lambda + 2 \, \mu > 0$: \\
Since we have
\begin{equation}
\langle \H^{\ZJ}(\sigma).D,D \rangle =
\langle\H^{\ZJ}(\mathcal{F}^{-1}(\sigma)).D,D \rangle
= \langle\H^{\ZJ}(B).D,D \rangle \qquad
\text{for} \qquad B = \mathcal{F}^{-1}(\sigma) \in \Sym^{++}(3),
\end{equation}
it follows (cf. \eqref{eqintro004})
%\blue{1. Zeile bitte ausführlicher: Bei der Herleitung ist mir die 1. Zeile unklar. Irgendwie erhalte ich andere Faktoren. Vermutlich aber mein Fehler!} \red{nochmal prüfen}
\begin{equation}
\label{eqthelatter01}
\begin{alignedat}{2}
\langle \H^{\ZJ}(B) . D , D \rangle &= \big\langle \frac{\mu}{2} \,
(D \, B + B \, D + B^{-1} (D \, B + B \, D) \, B^{-1}) , D \big \rangle + \frac{\lambda}{2} \,
\langle B^{-1} , D \, B + B \, D \rangle \cdot \big \langle \id, D \big\rangle \\
&= \frac{\mu}{2} \, (2 \, \langle B \, D , D \rangle + 2 \, \langle B^{-1}\, D,
D \rangle) + \frac{\lambda}{2} \, \langle 2\tr(D) \cdot \id , D \rangle \\
&= \mu \, ( \langle B \, D , D \rangle + \langle B^{-1} \, D, D \rangle ) +
\lambda \, \tr^2(D) = \mu \, \langle (B + B^{-1}) \, D, D \rangle + \lambda \,
\tr^2(D) \\
&\ge \mu \underbrace{\lambda_{\min}(B + B^{-1})}_{\ge 2, \;
\textnormal{see footnote\footnotemark}} \cdot \norm{D}^2 + \lambda \,
\tr^2(D) \ge \underbrace{2\mu \, \norm{D}^2 + \lambda \,
\tr^2(D)}_{> 0 \; \textnormal{for} \; \mu > 0, \, 2\mu + 3 \lambda > 0}
\ge c^+(\mu, \lambda) \cdot \norm{D}^2 \, .
\end{alignedat}
\end{equation}
\footnotetext
{
It is
\begin{equation} \begin{alignedat}{2}
\langle \xi, (B + B^{-1}) \, \xi \rangle &= \langle \xi,
(Q^T \, \text{diag} \, Q + Q^T \, \text{diag}^{-1} \, Q) \xi \rangle =
\langle Q \, \xi, (\text{diag} + \text{diag}^{-1}) \, Q \, \xi \rangle \\
&= \left\langle \eta, \text{diag}\left(\lambda_1 + \frac{1}{\lambda_1},
\lambda_2 + \frac{1}{\lambda_2}, \lambda_3 + \frac{1}{\lambda_3}\right) \,
\eta \right\rangle = \sum_{i=1}^3 \eta_i^2 \, \underbrace{
\left(\lambda_i + \frac{1}{\lambda_i}\right)}_{\ge 2} \ge 2 \,
\norm{\eta}^2 = 2 \, \norm{\xi}^2 \, .
\end{alignedat} \end{equation}
}

\subsubsection{Minor and major symmetry}\label{sec:minmajsym}
Minor symmetry of the induced fourth-order tangent stiffness tensor $\H^{\ZJ}(\sigma)$, i.e.
	\begin{align}
	D \in \Sym(3) \mapsto \H^{\ZJ}(\sigma). D \in \Sym(3),
	\end{align}
follows directly by using that
$B = \mathcal{F}^{-1}(\sigma) \in \Sym^{++}(3), \, B^{-1}$ and $D$ are
symmetric. \\
\\
For major symmetry of $\H^{\ZJ}(\sigma)$, i.e.
\begin{align}
\langle \H^{\ZJ}(\sigma). D_1 , D_2 \rangle = \langle \H^{\ZJ}(\sigma) . D_2,
D_1 \rangle \, ,
\end{align}
we calculate
\begin{equation}
\begin{alignedat}{2}
\langle \H^{\ZJ}(\sigma). D_1, D_2 \rangle &= \frac{\mu}{2} \, ( \langle B \, D_1 ,
D_2 \rangle + \langle D_1 \, B , D_2 \rangle + \langle D_1 \, B^{-1},
D_2 \rangle + \langle B^{-1} \, D_1 , D_2 \rangle ) + \lambda \, \tr(D_1)
\cdot \tr(D_2) \\
&= \frac{\mu}{2} \, ( \langle D_1, B \, D_2 \rangle + \langle D_1, D_2 \, B \rangle
+ \langle D_1, D_2 \, B^{-1} \rangle + \langle D_1, B^{-1} \, D_2 \rangle )
+ \lambda \, \tr(D_2) \cdot \tr(D_1) \\
&= \frac{\mu}{2} \, ( \langle B \, D_2 + D_2 \, B + B^{-1} \, D_2 + D_2 \, B^{-1}, D_1 \rangle )
+ \lambda \, \tr(D_2) \cdot \tr(D_1) = \langle \H^{\ZJ}(\sigma) . D_2,
D_1 \rangle.
\end{alignedat}
\end{equation}
\begin{Example}  Consider
        \begin{equation}
            \sigma(B) = \frac{B}{(\det B)^{\frac13}} - \id = (\det B)^{-\frac13} \cdot B - \id \, , \quad \sigma(\id) = 0 \, .
        \end{equation}
        Hence
	\begin{equation}
		\begin{alignedat}{2}
			\DD_B \sigma(B) . H &= -\frac13 \, (\det B)^{-\frac43} \cdot \langle \Cof B , H \rangle \cdot B + (\det B)^{-\frac13} \cdot H \\
			&= -\frac13 \, (\det B)^{-\frac43} \cdot \det B \cdot \langle B^{-1} , H \rangle \cdot B + (\det B)^{-\frac13} \cdot H \\
			&= -\frac13 \, (\det B)^{-\frac13} \cdot \langle B^{-1} , H \rangle \cdot B + (\det B)^{-\frac13} \cdot H \\
			\implies \quad \H^{\ZJ}(\sigma).D :&= \DD_B\sigma(B).[B \, D + D \, B] \\
			&= -\frac13 \, (\det B)^{-\frac13} \cdot \langle B^{-1}, B \, D + D \, B \rangle \, B + (\det B)^{-\frac13} \, (B \, D + D \, B) \\
			&= -\frac13 \, (\det B)^{-\frac13} \cdot 2 \, \tr(D) \cdot B + (\det B)^{-\frac13} \, (B \, D + D \, B) \\
			\implies \quad \langle \H^{\ZJ}(\sigma). D_1, D_2 \rangle &= -\frac13 \, (\det B)^{-\frac13} \cdot 2 \, \underbrace{\tr(D_1) \, \langle B , D_2 \rangle}_{\text{not interchangeable}} + (\det B)^{-\frac13} \cdot \underbrace{\langle B \, D_1 + D_1 \, B, D_2 \rangle}_{\text{interchangeable}} \, .
		\end{alignedat}
	\end{equation}
	So, this example shows that $\H^{\ZJ}(\sigma)$ is in general not major symmetric (not self-adjoint).
\end{Example}
\subsubsection{Monotonicity of $\widehat \sigma(\log B)$ in $\log B$}
\label{sec:app_monotonicity_log}
It is proven in \cite{CSP2024}, that even for a larger, more general class of corotational derivatives, the positive definiteness of a minor symmetric fourth-order tangent stiffness tensor $\H^{\ZJ}(\sigma)$ implies the strict monotonicity of $\widehat \sigma(\log B)$ in $\log B$, i.e. for all $B_1, B_2 \in \Sym^{++}(3)$ with $B_1 \neq B_2$ we have
	\begin{align}
	\sym \, \H^{\ZJ}(\sigma) \in \Sym^{++}_4(6) \qquad \implies \qquad \langle \widehat \sigma(\log B_1) - \widehat \sigma(\log B_2), \log B_1 - \log B_2 \rangle > 0.
	\end{align}
However, for the reader's convenience, we will prove strict monotonicity of $\widehat \sigma(\log B)$ in $\log B$ by direct inspection for our constitutive law. For simplicity, let us
only consider the case $\mu, \, \lambda > 0$ (auxetic-like response is excluded, since we assumed Poisson ratio $\nu > 0$). \\
%\blue{meaning? Was bedeutet auxetic?}\\
%
The term $B - \id$ is monotone in $\log B$, because
\begin{equation}
\label{eqmonotonesig}
\begin{alignedat}{2}
\langle \widehat{\sigma}(\log B_1) - \widehat{\sigma}(\log B_2),
\log B_1 - \log B_2 \rangle &= \langle \sigma(B_1) - \sigma(B_2) ,
\log B_1 - \log B_2 \rangle \\
&= \mu \, \langle B_1 - \id - (B_2 - \id), \log B_1 - \log B_2 \rangle \\
&= \mu \, \langle B_1 - B_2, \log B_1 - \log B_2 \rangle \\
&= \mu \, \langle \log B_1 - \log B_2, B_1 - B_2 \rangle > 0,
\end{alignedat}
\end{equation}
since $B \mapsto \log B$ is a monotone matrix function (cf. the upcoming paper
by Martin et al.~\cite{MartinVossGhibaNeff}). Furthermore, $\tr(\log B) \cdot \id$ is monotone in
$\log B$, since it is linear. \\
It remains to check the term $\sigma = \frac{1}{\mu} \sigma_2 = \id - B^{-1}$.
We write
\begin{equation}
\begin{alignedat}{2}
\langle \sigma(B_1) - \sigma(B_2), \log B_1 - \log B_2 \rangle &= \langle \id
- B_1^{-1} - (\id - B_2^{-1}), \log B_1 - \log B_2 \rangle \\
&= \langle -B_1^{-1} - (-B_2^{-1}), \log B_1 - \log B_2 \rangle \\
&= \langle -B_1^{-1} + B_2^{-1} , - \log B_1^{-1} - (-\log B_2^{-1}) \rangle \\
&= \langle B_2^{-1} - B_1^{-1} , \log B_2^{-1} - \log B_1^{-1} \rangle > 0 \\
&= \langle \log X - \log Y, X - Y \rangle > 0, \quad X=B_2^{-1},
\quad Y = B_1^{-1},
\end{alignedat}
\end{equation}
since $\log(\cdot)$ is monotone in its argument.

\subsection{Further observations for the constitutive law and the rate-formulation}
\subsubsection{Monotonicity in one dimension}
\begin{figure}[t]
\begin{center}
\begin{minipage}[h!]{0.95\linewidth}
\if\Bilder y
\centering
% vorher scale=0.25
\includegraphics[scale=0.2]{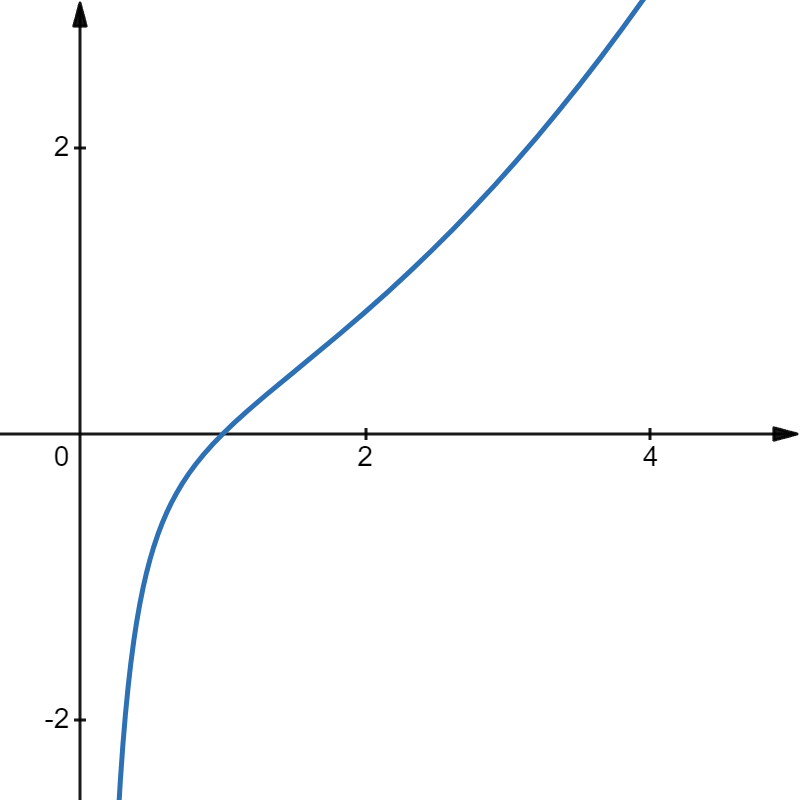}
\put(-7,60){\footnotesize{$\lambda$}}
\put(-168,147){\footnotesize{$\sigma(\lambda)$}}
\qquad
\includegraphics[scale=0.2]{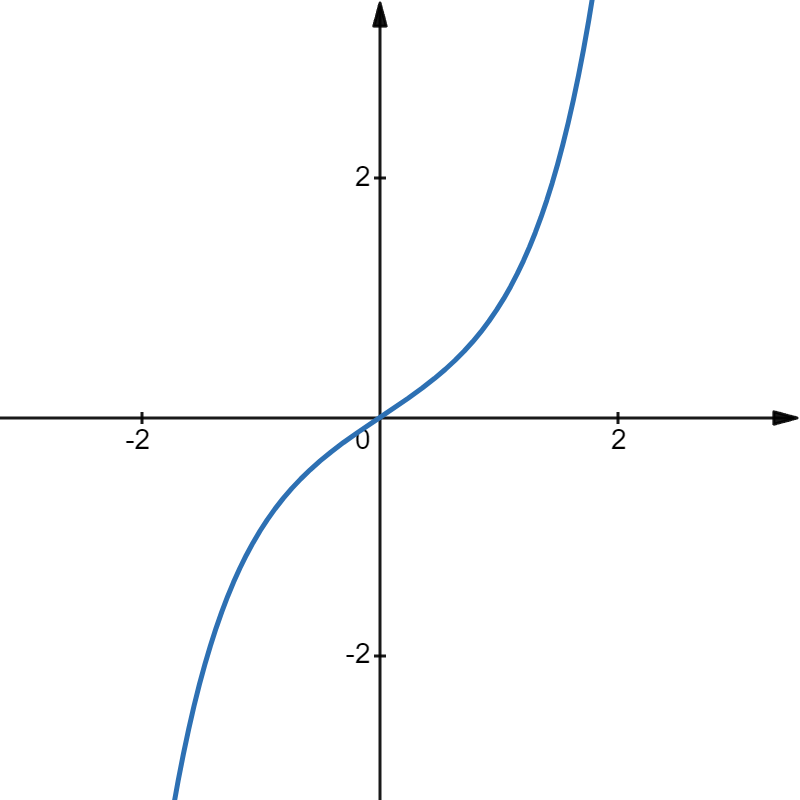}
\put(-7,63){\footnotesize{$\log \lambda$}}
\put(-120,147){\footnotesize{$\widehat{\sigma}(\log \lambda)$}}
\fi
\caption{We consider $\sigma: \R^+ \to \R$ and set
$\lambda \mapsto \sigma(\lambda) := \widehat{\sigma}(\log \lambda)$. Depicted
is the graph of the one-dimensional forces
$\sigma(\lambda)=\frac16\left(\lambda^2-
\frac{1}{\lambda^2}+\log\lambda^2 \right)$ (left) and
$\widehat{\sigma}(\log \lambda) = \frac16 \, ( 2 \, \sinh(2 \, \log \lambda) + 2 \,
\log \lambda)$ (right), modelled after our Neo-Hooke type law
$\sigma(B)=\frac{\mu}{2}\, (B-B^{-1})+\frac{\lambda}{2}\,
(\log\det B) \cdot\id$. Here, $\lambda \mapsto \sigma(\lambda)$ and \break
$\log \lambda \mapsto \widehat{\sigma}(\log \lambda)$ are both monotone.
In the multi-dimensional setting, $B \mapsto \sigma(B)$ is not monotone, while
$\log B \mapsto \widehat{\sigma}(\log B)$ is monotone.}
\label{fig3}
\end{minipage}
\end{center}
\end{figure}

\begin{figure}[t]
\begin{center}
\begin{minipage}[h!]{0.8\linewidth}
\if\Bilder y
\centering
\includegraphics[scale=0.3]{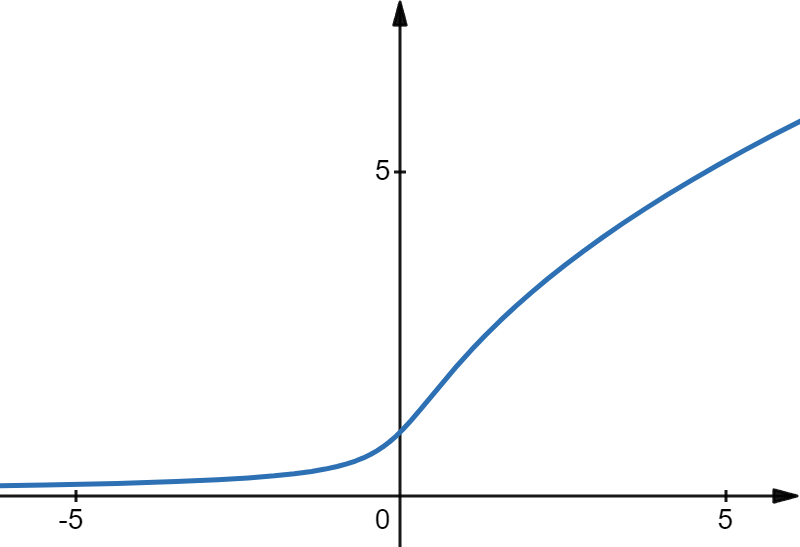}
\put(-5,2){$\sigma$}
\put(-150,110){$\lambda = \mathcal{F}^{-1}(\sigma)$}
\fi
\caption{It can be seen that
$\mathcal{F}^{-1}: \R \to \R^+, \; \lambda = \mathcal{F}^{-1}(\sigma)$ is the
smooth inverse mapping. The same is true for
$\mathcal{F}^{-1}: \Sym(3) \to \Sym^{++}(3)$ for the constitutive law
$\sigma(B) = \frac{\mu}{2} \, (B-B^{-1}) + \frac{\lambda}{2} \,
(\log \det B) \cdot \id$.}
\label{fig9}
\end{minipage}
\end{center}
\end{figure}

As we have shown in Section \ref{sec:app_monotonicity_log},
%a previous section of the Appendix \ref{secA1}, 
the constitutive law
\begin{align}
\sigma: \Sym^{++}(3) \to \Sym(3), \qquad B \mapsto \sigma(B) = \frac{\mu}{2} \,
(B - B^{-1}) + \frac{\lambda}{2} \, \log \det B \cdot \id
\end{align}
fulfills monotonicity of $\log B \mapsto \widehat \sigma(\log B) = \sigma(B)$ in $\log B$. Note carefully that in three dimensions this is not equivalent to monotonicity of $B \mapsto \sigma(B)$ in $B$, as was discussed in Section \ref{isocauchyelast} (cf.~\cite{CSP2024}). However, in one dimension, both types of monotonicity are equivalent, as illustrated by the Figures \ref{fig3} and \ref{fig9}. \\

%
%\clearpage
%
\subsubsection{Geometrical meaning of the Eulerian rate of deformation $D$}
As a recapitulation for the reader we provide the following identities
\begin{equation} \begin{alignedat}{2}
\delta(t) &= \norm{\varphi(x_1,t) - \varphi(x_2,t)}, \\
\frac12 \delta^2(t) &= \frac12 \langle \varphi(x_1,t) - \varphi(x_2,t),
\varphi(x_1,t) - \varphi(x_2,t) \rangle, \\
\frac{\dif}{\dif t} \, \frac12 \delta^2(t) &= \frac{\dif}{\dif t} \, \frac12
\langle \varphi(x_1,t) - \varphi(x_2,t), \varphi(x_1,t) - \varphi(x_2,t)
\rangle = \langle \dot{\varphi}(x_1,t) - \dot{\varphi}(x_2,t),
\varphi(x_1,t) - \varphi(x_2,t) \rangle, \\
\delta(t) \cdot \dot{\delta}(t) &= \langle v(\xi_1,t) - v(\xi_2,t),
\xi_1 - \xi_2 \rangle \approx \langle \DD_{\xi} v(\xi_1) . \delta \xi,
\delta \xi \rangle, \quad \text{if} \quad \norm{\xi_1 - \xi_2} \ll 1,
\qquad \xi_2 = \xi_1 + \delta \xi, \\
\dot{\delta}(t) &\cong \frac{1}{\delta(t)} \cdot \langle
[\sym \, \DD_{\xi} v(\xi)] . \delta \xi, \delta \xi \rangle, \qquad
D = \sym \, \DD_{\xi} v(\xi) = \mathbb{S}(\sigma) .
\frac{\DD^{\ZJ}}{\DD t}[\sigma(\xi,t)] \, .
\end{alignedat} \end{equation}
Thus we observe that $D$ describes the change of stretch per unit length, hence the
name ``stretching'' for \break $D = \sym \, \DD_{\xi} v(\xi)$. Here,
$\mathbb{S}(\sigma) = [\H^{\ZJ}(\sigma)]^{-1}$ is the induced fourth-order compliance tensor in the rate-formulation.

\end{appendix}

\end{document}